%
%
%
\newif\ifhypers                 
\hypersfalse
\ifhypers

\fi
%

\ifx\begin\undefined\else\global\advance\srcdepth by
1\expandafter \fi

\def\begin{}
\newcount\srcdepth
\srcdepth=1

\outer\def\bye{\global\advance\srcdepth by -1
  \ifnum\srcdepth=0
    \def\endcmd{\vfill\supereject\nopagenumbers\par\vfill\supereject\end}
  \else\def\endcmd{}\fi
  \endcmd
}



\def\initialize#1#2#3#4#5#6{
  \ifnum\srcdepth=1
  \magnification=#1
  \hsize = #2
  \vsize = #3
  \hoffset=#4
  \advance\hoffset by -\hsize
  \divide\hoffset by 2
  \advance\hoffset by -1truein
  \voffset=#5
  \advance\voffset by -\vsize
  \divide\voffset by 2
  \advance\voffset by -1truein
  \advance\voffset by #6
  \baselineskip=13pt
  \emergencystretch = 0.05\hsize
  \fi
}

\def\print{\initialize{1095}
  {5.5truein}{8.5truein}{8.5truein}{11truein}{-.0625truein}}

\newif\ifblackboardbold

\blackboardboldtrue


\font\titlefont=cmbx12 scaled\magstephalf
\font\sectionfont=cmbx12

\font\scriptit=cmti10 at 7pt
\font\scriptsl=cmsl10 at 7pt
\scriptfont\itfam=\scriptit
\scriptfont\slfam=\scriptsl


\newfam\bboldfam
\ifblackboardbold
\font\tenbbold=msbm10
\font\sevenbbold=msbm7
\font\fivebbold=msbm5
\textfont\bboldfam=\tenbbold
\scriptfont\bboldfam=\sevenbbold
\scriptscriptfont\bboldfam=\fivebbold
\def\bbold{\fam\bboldfam\tenbbold}
\else
\def\bbold{\bf}
\fi


\newfam\msamfam
\font\tenmsam=msam10
\font\sevenmsam=msam7
\font\fivemsam=msam5
\textfont\msamfam=\tenmsam
\scriptfont\msamfam=\sevenmsam
\scriptscriptfont\msamfam=\fivemsam

\newfam\msbmfam
\font\tenmsbm=msam10
\font\sevenmsbm=msam7
\font\fivemsbm=msam5
\textfont\msbmfam=\tenmsbm
\scriptfont\msbmfam=\sevenmsbm
\scriptscriptfont\msbmfam=\fivemsbm

\newcount\amsfamcount 
\newcount\classcount   
\newcount\positioncount
\newcount\codecount
\newcount\n             
\def\newsymbol#1#2#3#4#5{               
\n="#2                                  
\ifnum\n=1 \amsfamcount=\msamfam\else   
\ifnum\n=2 \amsfamcount=\msbmfam\else   
\ifnum\n=3 \amsfamcount=\eufmfam
\fi\fi\fi
\multiply\amsfamcount by "100           
\classcount="#3                 
\multiply\classcount by "1000           
\positioncount="#4#5            
\codecount=\classcount                  
\advance\codecount by \amsfamcount      
\advance\codecount by \positioncount
\mathchardef#1=\codecount}              


\font\Arm=cmr9
\font\Ai=cmmi9
\font\Asy=cmsy9
\font\Abf=cmbx9
\font\Brm=cmr7
\font\Bi=cmmi7
\font\Bsy=cmsy7
\font\Bbf=cmbx7
\font\Crm=cmr6
\font\Ci=cmmi6
\font\Csy=cmsy6
\font\Cbf=cmbx6

\ifblackboardbold
\font\Abbold=msbm10 at 9pt
\font\Bbbold=msbm7
\font\Cbbold=msbm5 at 6pt
\fi

\def\small{%
\textfont0=\Arm \scriptfont0=\Brm \scriptscriptfont0=\Crm
\textfont1=\Ai \scriptfont1=\Bi \scriptscriptfont1=\Ci
\textfont2=\Asy \scriptfont2=\Bsy \scriptscriptfont2=\Csy
\textfont\bffam=\Abf \scriptfont\bffam=\Bbf \scriptscriptfont\bffam=\Cbf
\def\rm{\fam0\Arm}\def\mit{\fam1}\def\oldstyle{\fam1\Ai}%
\def\bf{\fam\bffam\Abf}%
\ifblackboardbold
\textfont\bboldfam=\Abbold
\scriptfont\bboldfam=\Bbbold
\scriptscriptfont\bboldfam=\Cbbold
\def\bbold{\fam\bboldfam\Abbold}%
\fi
\rm
}








\newlinechar=`@
\def\forwardmsg#1#2#3{\immediate\write16{@*!*!*!* forward reference should
be: @\noexpand\forward{#1}{#2}{#3}@}}
\def\nodefmsg#1{\immediate\write16{@*!*!*!* #1 is an undefined reference@}}

\def\forwardsub#1#2{\def\newref{{#2}{#1}}}

\def\forward#1#2#3{%
\expandafter\expandafter\expandafter\forwardsub\expandafter{#3}{#2}
\expandafter\ifx\csname#1\endcsname\relax\else%
\expandafter\ifx\csname#1\endcsname\newref\else%
\forwardmsg{#1}{#2}{#3}\fi\fi%
\expandafter\let\csname#1\endcsname\newref}

\def\firstarg#1{\expandafter\argone #1}\def\argone#1#2{#1}
\def\secondarg#1{\expandafter\argtwo #1}\def\argtwo#1#2{#2}

\def\ref#1{\expandafter\ifx\csname#1\endcsname\relax
  {\nodefmsg{#1}\bf`#1'}\else
  \expandafter\firstarg\csname#1\endcsname
  ~\htmllocref{#1}{\expandafter\secondarg\csname#1\endcsname}\fi}

\def\refs#1{\expandafter\ifx\csname#1\endcsname\relax
  {\nodefmsg{#1}\bf`#1'}\else
  \expandafter\firstarg\csname #1\endcsname
  s~\htmllocref{#1}{\expandafter\secondarg\csname#1\endcsname}\fi}

\def\refn#1{\expandafter\ifx\csname#1\endcsname\relax
  {\nodefmsg{#1}\bf`#1'}\else
  \htmllocref{#1}{\expandafter\secondarg\csname #1\endcsname}\fi}



\def\widow#1{\vskip 0pt plus#1\vsize\goodbreak\vskip 0pt plus-#1\vsize}

%
\def\begincolor#1{\ifhypers\fi}%
\def\endcolor{\ifhypers\fi}%
\def\colored#1#2{%
  \begincolor{#1}#2\endcolor}
\def\linkcolor{[1.0 0.1 0.1]}%
\def\sectioncolor{[0 0.4 0.8]}%
\def\proclaimcolor{[0 0.4 0.8]}%
\def\proofcolor{[0 0 0]}%
%
%
%
%
\edef\freehash{\catcode`\noexpand\#=\the\catcode`\#}%
\catcode`\#=12
\freehash
\let\freehash=\relax
\def\puthtml#1{\ifhypers\fi}
\def\htmlanchor#1#2{\puthtml{<a name="#1">}#2\puthtml{</a>}}

\def\@pdfm@mark#1{}
\def\setlink#1{\colored{\linkcolor}{#1}}%
%
%
%
%
\def\htmllocref#1#2{\ifhypers\leavevmode\fi\setlink{#2}\ifhypers\fi\relax}%
%
%
%
%
\def\Acrobatmenu#1#2{%
  \@pdfm@mark{%
    bann <<
      /Type /Annot
      /Subtype /Link
      /A <<
        /S /Named
        /N /#1
      >>
      /Border [\@pdfborder]
      /C [\@menubordercolor]
    >>%
   }%
  \Hy@colorlink{\@menucolor}#2\Hy@endcolorlink
  \@pdfm@mark{eann}%
}
\def\@pdfborder{0 0 1}
\def\@menubordercolor{1 0 0}
\def\@menucolor{red}
%
%



\def\marginlabel#1{}

\def\showlabelsabove{
\font\labelfont=cmss10 at 6pt
\def\marginlabel##1{\rlap{\smash{\raise 10pt\hbox{\labelfont##1}}}}
}

\newcount\seccount
\newcount\proccount
\seccount=0
\proccount=0

\def\stdskip{\vskip 9pt plus3pt minus 3pt}
\def\stdbreak{\par\removelastskip\penalty-100\stdskip}

\def\proof{\stdbreak\noindent \colored{\proofcolor}{{\sl Proof. }}}

\def\qed{\vrule height 1.2ex width .9ex depth .1ex}

\def\Box{
  \ifmmode\eqno\qed
  \else\ifvmode\removelastskip\line{\hfil\qed}
  \else\unskip\quad\hskip-\hsize
    \hbox{}\hskip\hsize minus 1em\qed\par
  \fi\stdbreak\fi}

\def\references{
  \removelastskip
  \widow{.05}
  \vskip 24pt plus 6pt minus 6 pt
  \parindent=0pt
  \frenchspacing
  \leftline{\sectionfont \colored{\sectioncolor}{References}}
  \ifhypers\global\advance\seccount by 1\global\proccount=0\relax\edef\numtoks{\number\seccount}%
  \hbox{}\fi
  \nobreak\stdskip\noindent}

\def\ifempty#1#2\endB{\ifx#1\endA}
\def\makeref#1#2#3{\ifempty#1\endA\endB\else\forward{#1}{#2}{#3}\fi}

\outer\def\section#1 #2\par{
  \removelastskip
  \global\advance\seccount by 1
  \global\proccount=0\relax
                \edef\numtoks{\number\seccount}
  \htmlanchor{#1}{\makeref{#1}{Section}{\numtoks}}
  \widow{.05}
  \vskip 24pt plus 6pt minus 6 pt
  \message{#2}
  \ifhypers\hbox{}\fi
  \leftline{\marginlabel{#1}\sectionfont\colored{\sectioncolor}{\numtoks}\quad \colored{\sectioncolor}{#2}}
  \nobreak\stdskip}

\def\proclamation#1#2{
  \outer\expandafter\def\csname#1\endcsname##1 ##2\par{
  \stdbreak
  \global\advance\proccount by 1
  \edef\numtoks{\number\seccount.\number\proccount}
  \htmlanchor{##1}{\makeref{##1}{#2}{\numtoks}}
  \noindent{\marginlabel{##1}\bf \colored{\proclaimcolor}{#2} \colored{\proclaimcolor}{\numtoks}\enspace}
  {\sl##2\par}
  \stdbreak}}

\def\othernumbered#1#2{
  \outer\expandafter\def\csname#1\endcsname##1{
  \stdbreak
  \global\advance\proccount by 1
  \edef\numtoks{\number\seccount.\number\proccount}
  \htmlanchor{##1}{\makeref{##1}{#2}{\numtoks}}
  \noindent{\marginlabel{##1}\bf \colored{\proclaimcolor}{#2} \colored{\proclaimcolor}{\numtoks}\enspace}}}

\proclamation{definition}{Definition}
\proclamation{lemma}{Lemma}
\proclamation{proposition}{Proposition}
\proclamation{theorem}{Theorem}
\proclamation{corollary}{Corollary}
\proclamation{conjecture}{Conjecture}

\othernumbered{example}{Example}
\othernumbered{remark}{Remark}
\othernumbered{construction}{Construction}
\othernumbered{problem}{Problem}

\def\figure#1{
 \global\advance\figcount by 1
 \goodbreak
 \midinsert#1\smallskip
 \centerline{Figure~\number\figcount}
 \endinsert}

\def\capfigure#1#2{
 \global\advance\figcount by 1
 \goodbreak
 \midinsert#1\smallskip
 \vbox{\small\noindent {\bf Figure~\number\figcount:} #2}
 \endinsert}

\def\capfigurepair#1#2#3#4{
 \goodbreak
 \midinsert
 #1\smallskip
 \global\advance\figcount by 1
 \vbox{\small\noindent {\bf Figure~\number\figcount:} #2}
 \vskip 12pt
 #3\smallskip
 \global\advance\figcount by 1
 \vbox{\small\noindent {\bf Figure~\number\figcount:} #4}
 \endinsert}


\def\baretable#1#2{
\vbox{\offinterlineskip\halign{
 \strut\kern #1\hfil##\kern #1
 &&\kern #1\hfil##\kern #1\cr
 #2
}}}

\def\gridtablesub#1#2#3{
\vbox{\offinterlineskip\halign{
 \strut\vrule\kern #1\hfil##\hfil\kern #2\vrule
 &&\kern #1\hfil##\kern #2\vrule\cr
 \noalign{\hrule}
 #3
 \noalign{\hrule}
}}}






\newif\iftextures
\input epsf

\newcount\figcount
\figcount=0
\newcount\figxscale
\newcount\figyscale
\newcount\figxoffset
\newcount\figyoffset
\newbox\drawing
\newcount\drawbp
\newdimen\drawx
\newdimen\drawy
\newdimen\ngap
\newdimen\sgap
\newdimen\wgap
\newdimen\egap

\def\drawbox#1#2#3{\vbox{
  \epsfgetbb{#2.eps} 
  \drawbp=\epsfurx
  \advance\drawbp by-\epsfllx\relax
  \multiply\drawbp by #1
  \divide\drawbp by 100
  \drawx=\drawbp bp
  \drawbp=\epsfury
  \advance\drawbp by-\epsflly\relax
  \multiply\drawbp by #1
  \divide\drawbp by 100
  \drawy=\drawbp bp
  \iftextures
  		\figxscale=#1
    \multiply\figxscale by 10
    \setbox\drawing=\vbox to \drawy{\vfil
      \hbox to \drawx{\special{illustration #2.eps scaled
\number\figxscale}\hfil}}
  \else 
    \figxoffset=-\epsfllx
    \multiply\figxoffset by#1
    \divide\figxoffset by100
    \figyoffset=-\epsflly
    \multiply\figyoffset by#1
    \divide\figyoffset by100
    \setbox\drawing=\vbox to \drawy{\vfil
      \hbox to \drawx{\includegraphics{#2.eps}\hfil}}
  \fi
  \setbox\drawing=\vbox{\offinterlineskip\box\drawing\kern 0pt}
   \ngap=0pt \sgap=0pt \wgap=0pt \egap=0pt
  \setbox0=\vbox{\offinterlineskip
    \box\drawing \ifgridlines\drawgrid\drawx\drawy\fi #3}
  \kern\ngap\hbox{\kern\wgap\box0\kern\egap}\kern\sgap}}

\def\draw#1#2#3{
  \setbox\drawing=\drawbox{#1}{#2}{#3}
  \global\advance\figcount by 1
  \edef\numtoks{\number\figcount}
  \makeref{fig:#2}{Figure}{\numtoks}
  \goodbreak
  \midinsert
  \centerline{\ifgridlines\boxgrid\drawing\fi\box\drawing}
  \smallskip
  \vbox{\offinterlineskip
    \centerline{Figure~\number\figcount}
    \smash{\marginlabel{#2}}}
  \endinsert}

\def\capdraw#1#2#3#4{
  \setbox\drawing=\drawbox{#1}{#2}{#3}
  \global\advance\figcount by 1
  \edef\numtoks{\number\figcount}
  \makeref{fig:#2}{Figure}{\numtoks}
  \goodbreak
  \midinsert
  \centerline{\ifgridlines\boxgrid\drawing\fi\box\drawing}
  \smallskip
  \vbox{\offinterlineskip
    \vskip 4pt
    \vbox{\centerline{\lineskip=3pt\small\noindent
                       {\bf Figure~\number\figcount:} #4}}
    \smash{\marginlabel{fig:#2}}}
  \endinsert}

\def\capdrawpair#1#2#3#4#5#6#7#8{
  \goodbreak
  \midinsert
  \setbox\drawing=\drawbox{#1}{#2}{#3}
  \global\advance\figcount by 1
  \edef\numtoks{\number\figcount}
  \makeref{fig:#2}{Figure}{\numtoks}
  \centerline{\ifgridlines\boxgrid\drawing\fi\box\drawing}
  \smallskip
  \vbox{\offinterlineskip
    \vskip 4pt
    \vbox{\lineskip=3pt\small\noindent {\bf Figure~\number\figcount:} #4}
    \smash{\marginlabel{fig:#2}}}
  \vskip 12pt
  \setbox\drawing=\drawbox{#5}{#6}{#7}
  \global\advance\figcount by 1
  \edef\numtoks{\number\figcount}
  \makeref{fig:#6}{Figure}{\numtoks}
  \centerline{\ifgridlines\boxgrid\drawing\fi\box\drawing}
  \smallskip
  \vbox{\offinterlineskip
    \vskip 4pt
    \vbox{\lineskip=3pt\small\noindent {\bf Figure~\number\figcount:} #8}
    \smash{\marginlabel{fig:#6}}}
  \endinsert}

\def\nextfigtoks{%
  \advance\figcount by 1%
  \edef\numtoks{\number\figcount}%
  \advance\figcount by -1}

\newif\ifgridlines
\newbox\figtbox
\newbox\figgbox
\newdimen\figtx
\newdimen\figty

\newdimen\bwd
\bwd=2sp 

\def\hline#1{\vbox{\smash{\hbox to #1{\leaders\hrule height \bwd\hfil}}}}

\def\vline#1{\hbox to 0pt{%
  \hss\vbox to #1{\leaders\vrule width \bwd\vfil}\hss}}

\def\clap#1{\hbox to 0pt{\hss#1\hss}}
\def\vclap#1{\vbox to 0pt{\offinterlineskip\vss#1\vss}}

\def\hstutter#1#2{\hbox{%
  \setbox0=\hbox{#1}%
  \hbox to #2\wd0{\leaders\box0\hfil}}}

\def\vstutter#1#2{\vbox{
  \setbox0=\vbox{\offinterlineskip #1}
  \dp0=0pt
  \vbox to #2\ht0{\leaders\box0\vfil}}}

\def\crosshairs#1#2{
  \dimen1=.002\drawx
  \dimen2=.002\drawy
  \ifdim\dimen1<\dimen2\dimen3\dimen1\else\dimen3\dimen2\fi
  \setbox1=\vclap{\vline{2\dimen3}}
  \setbox2=\clap{\hline{2\dimen3}}
  \setbox3=\hstutter{\kern\dimen1\box1}{4}
  \setbox4=\vstutter{\kern\dimen2\box2}{4}
  \setbox1=\vclap{\vline{4\dimen3}}
  \setbox2=\clap{\hline{4\dimen3}}
  \setbox5=\clap{\copy1\hstutter{\box3\kern\dimen1\box1}{6}}
  \setbox6=\vclap{\copy2\vstutter{\box4\kern\dimen2\box2}{6}}
  \setbox1=\vbox{\offinterlineskip\box5\box6}
  \smash{\vbox to #2{\hbox to #1{\hss\box1}\vss}}}

\def\boxgrid#1{\rlap{\vbox{\offinterlineskip
  \setbox0=\hline{\wd#1}
  \setbox1=\vline{\ht#1}
  \smash{\vbox to \ht#1{\offinterlineskip\copy0\vfil\box0}}
  \smash{\vbox{\hbox to \wd#1{\copy1\hfil\box1}}}}}}

\def\drawgrid#1#2{\vbox{\offinterlineskip
  \dimen0=\drawx
  \dimen1=\drawy
  \divide\dimen0 by 10
  \divide\dimen1 by 10
  \setbox0=\hline\drawx
  \setbox1=\vline\drawy
  \smash{\vbox{\offinterlineskip
    \copy0\vstutter{\kern\dimen1\box0}{10}}}
  \smash{\hbox{\copy1\hstutter{\kern\dimen0\box1}{10}}}}}

\def\figtext#1#2#3#4#5{
  \setbox\figtbox=\vbox{\hbox{#5}\kern 0pt}
  \figtx=-#3\wd\figtbox \figty=-#4\ht\figtbox
  \advance\figtx by #1\drawx \advance\figty by #2\drawy
  \dimen0=\figtx \advance\dimen0 by\wd\figtbox \advance\dimen0 by-\drawx
  \ifdim\dimen0>\egap\global\egap=\dimen0\fi
  \dimen0=\figty \advance\dimen0 by\ht\figtbox \advance\dimen0 by-\drawy
  \ifdim\dimen0>\ngap\global\ngap=\dimen0\fi
  \dimen0=-\figtx
  \ifdim\dimen0>\wgap\global\wgap=\dimen0\fi
  \dimen0=-\figty
  \ifdim\dimen0>\sgap\global\sgap=\dimen0\fi
  \smash{\rlap{\vbox{\offinterlineskip
    \hbox{\hbox to \figtx{}\ifgridlines\boxgrid\figtbox\fi\box\figtbox}
    \vbox to \figty{}
    \ifgridlines\crosshairs{#1\drawx}{#2\drawy}\fi
    \kern 0pt}}}}


\def\hpad#1#2#3{\hbox{\kern #1\hbox{#3}\kern #2}}
\def\vpad#1#2#3{\setbox0=\hbox{#3}\vbox{\kern #1\box0\kern #2}}




\def\stack#1#2#3{\vbox{\offinterlineskip
  \setbox2=\hbox{#2}
  \setbox3=\hbox{#3}
  \dimen0=\ifdim\wd2>\wd3\wd2\else\wd3\fi
  \hbox to \dimen0{\hss\box2\hss}
  \kern #1
  \hbox to \dimen0{\hss\box3\hss}}}


\def\hexp#1{%
  \setbox0=\hbox{${}^{#1}$}%
  \hbox to .5\wd0{\box0\hss}}

\def\hsub#1{%
  \setbox0=\hbox{${}_{#1}$}%
  \hbox to .5\wd0{\box0\hss}}



\def\bmatrix#1#2{{\left(\vcenter{\halign
  {&\kern#1\hfil$##\mathstrut$\kern#1\cr#2}}\right)}}

\def\rightarrowmat#1#2#3{
  \setbox1=\hbox{\small\kern#2$\bmatrix{#1}{#3}$\kern#2}
  \,\vbox{\offinterlineskip\hbox to\wd1{\hfil\copy1\hfil}
    \kern 3pt\hbox to\wd1{\rightarrowfill}}\,}

\def\leftarrowmat#1#2#3{
  \setbox1=\hbox{\small\kern#2$\bmatrix{#1}{#3}$\kern#2}
  \,\vbox{\offinterlineskip\hbox to\wd1{\hfil\copy1\hfil}
    \kern 3pt\hbox to\wd1{\leftarrowfill}}\,}

\def\rightarrowbox#1#2{
  \setbox1=\hbox{\kern#1\hbox{\small #2}\kern#1}
  \,\vbox{\offinterlineskip\hbox to\wd1{\hfil\copy1\hfil}
    \kern 3pt\hbox to\wd1{\rightarrowfill}}\,}

\def\leftarrowbox#1#2{
  \setbox1=\hbox{\kern#1\hbox{\small #2}\kern#1}
  \,\vbox{\offinterlineskip\hbox to\wd1{\hfil\copy1\hfil}
    \kern 3pt\hbox to\wd1{\leftarrowfill}}\,}








\def\quiremacro#1#2#3#4#5#6#7#8#9{
  \expandafter\def\csname#1\endcsname##1{
  \ifnum\srcdepth=1
  \magnification=#2
  \input quire
  \hsize=#3
  \vsize=#4
  \htotal=#5
  \vtotal=#6
  \shstaplewidth=#7
  \shstaplelength=#8
  \hoffset=\htotal
  \advance\hoffset by -\hsize
  \divide\hoffset by 2
  \ifnum\vsize<\vtotal
    \voffset=\vtotal
    \advance\voffset by -\vsize
    \divide\voffset by 2
  \fi
  \advance\voffset by #9
  \shhtotal=2\htotal
  \baselineskip=13pt
  \emergencystretch = 0.05\hsize
  \horigin=0.0truein
  \vorigin=0.0truein
  \shthickness=0pt
  \shoutline=0pt
  \shcrop=0pt
  \shvoffset=-1.0truein
  \ifnum##1>0\quire{#1}\else\qtwopages\fi
  \fi
}}



\quiremacro{letterbooklet} 
{1000}{4.79452truein}{7truein}{5.5truein}{8.5truein}{0.01pt}{0.66truein}
{-.0625truein}

\quiremacro{Afourbooklet}
{1095}{5.25truein}{7truein}{421truept}{595truept}{0.01pt}{0.66truein}
{-.0625truein}

\quiremacro{legalbooklet}
{1095}{5.25truein}{7truein}{7.0truein}{8.5truein}{0.01pt}{0.66truein}
{-.0625truein}

\quiremacro{twoupsub} 
{895}{4.5truein}{7truein}{5.5truein}{8.5truein}{0pt}{0pt}{.0625truein}


\quiremacro{Afourviewsub} 
{1000}{5.0228311in}{7.7625571in}{421truept}{595truept}{0.1pt}{0.5\vtotal}
{-.0625truein}


\quiremacro{viewsub}
{1095}{5.5truein}{8.5truein}{461truept}{666truept}{0.1pt}{0.5\vtotal}
{-.125truein}


\newcount\countA
\newcount\countB
\newcount\countC

\def\monthname{\begingroup
  \ifcase\number\month
    \or January\or February\or March\or April\or May\or June\or
    July\or August\or September\or October\or November\or December\fi
\endgroup}

\def\dayname{\begingroup
  \countA=\number\day
  \countB=\number\year
  \advance\countA by 0 
  \advance\countA by \ifcase\month\or
    0\or 31\or 59\or 90\or 120\or 151\or
    181\or 212\or 243\or 273\or 304\or 334\fi
  \advance\countB by -1995
  \multiply\countB by 365
  \advance\countA by \countB
  \countB=\countA
  \divide\countB by 7
  \multiply\countB by 7
  \advance\countA by -\countB
  \advance\countA by 1
  \ifcase\countA\or Sunday\or Monday\or Tuesday\or Wednesday\or
    Thursday\or Friday\or Saturday\fi
\endgroup}

\def\timename{\begingroup
   \countA = \time
   \divide\countA by 60
   \countB = \countA
   \countC = \time
   \multiply\countA by 60
   \advance\countC by -\countA
   \ifnum\countC<10\toks1={0}\else\toks1={}\fi
   \ifnum\countB<12 \toks0={\sevenrm AM}
     \else\toks0={\sevenrm PM}\advance\countB by -12\fi
   \relax\ifnum\countB=0\countB=12\fi
   \hbox{\the\countB:\the\toks1 \the\countC \thinspace \the\toks0}
\endgroup}

\def\timestamp{\dayname, \the\day\ \monthname\ \the\year, \timename}


\print



\def\COMMENT#1\par{\bigskip\hrule\smallskip#1\smallskip\hrule\bigskip}

\def\enma#1{{\ifmmode#1\else$#1$\fi}}

\def\mathbb#1{{\bbold #1}}
\def\mathbf#1{{\bf #1}}


\def\NN{\enma{\mathbb{N}}}

\def\PP{\enma{\mathbb{P}}}
\def\HH{\enma{\mathbb{H}}}
\def\ZZ{\enma{\mathbb{Z}}}


\def\cEE{\enma{\cal E}}
\def\cHH{\enma{\cal H}}
\def\cFF{\enma{\cal F}}
\def\cII{\enma{\cal I}}
\def\cKK{\enma{\cal K}}
\def\cOO{\enma{\cal O}}


\def\LL{\enma{\mathbf{L}}}

\def\boldN{\enma{\mathbf{N}}}
\font\tengoth=eufm10  \font\fivegoth=eufm5
\font\sevengoth=eufm7
\newfam\gothfam  \scriptscriptfont\gothfam=\fivegoth 
\textfont\gothfam=\tengoth \scriptfont\gothfam=\sevengoth


\def\codim{\mathop{\rm cod}\nolimits}

\def\lspan{{\mathop{\rm span}\nolimits}}

\def\ker{\mathop{\rm ker}\nolimits}
\def\deg{\mathop{\rm deg}\nolimits}

\def\length{{\mathop{\rm length}\nolimits}}
\def\rank{\mathop{\rm rank}\nolimits}
\def\link{\mathop{\rm link}\nolimits}
\def\rank{\mathop{\rm rank}\nolimits}
\def\dim{\mathop{\rm dim}\nolimits}

\def\Tor{\mathop{\rm Tor}\nolimits}

\def\Sym{\mathop{\rm Sym}\nolimits}

\def\Pic{\mathop{\rm Pic}\nolimits}
\def\cliff{{\mathop{\rm Cliff}\nolimits}}
\def\Gr{\enma{\rm Gr}}

\def\SL{\enma{\rm SL}}

\def\Sec{{\mathop{\rm Sec}\nolimits}}
%

%
%
\def\Eprime{\mathop{{}^\prime\!E}\nolimits}
\def\Edprime{\mathop{{}^{\prime\prime}\!E}\nolimits}

\newsymbol\boxtimes1202


\input diagrams.tex
\overfullrule=0pt
\proclamation{question}{Question}

\def\codim{\mathop{\rm cod}\nolimits}

\def\lspan{{\mathop{\rm span}\nolimits}}

\def\ker{\mathop{\rm ker}\nolimits}
\def\deg{\mathop{\rm deg}\nolimits}

\def\length{{\mathop{\rm length}\nolimits}}
\def\rank{\mathop{\rm rank}\nolimits}
\def\link{\mathop{\rm link}\nolimits}
\def\genus{\mathop{\rm genus}\nolimits}
\def\star{\mathop{\rm star}\nolimits}
\def\rank{\mathop{\rm rank}\nolimits}
\def\dim{\mathop{\rm dim}\nolimits}

\def\Tor{\mathop{\rm Tor}\nolimits}

\def\Sym{\mathop{\rm Sym}\nolimits}

\def\Pic{\mathop{\rm Pic}\nolimits}
\def\cliff{{\mathop{\rm Cliff}\nolimits}}
\def\Gr{\enma{\rm Gr}}

\def\SL{\enma{\rm SL}}
\def\Sec{{\mathop{\rm Sec}\nolimits}}

%

\def\Eprime{\mathop{{}^\prime\!E}\nolimits}
\def\Edprime{\mathop{{}^{\prime\prime}\!E}\nolimits}


\def\enma#1{{\ifmmode#1\else$#1$\fi}}
\def\mathbb#1{{\bbold #1}}
\def\mathbf#1{{\bf #1}}


\def\NN{\enma{\mathbb{N}}}

\def\PP{\enma{\mathbb{P}}}
\def\HH{\enma{\mathbf{H}}}

\def\ZZ{\enma{\mathbb{Z}}}


\def\cEE{\enma{\cal E}}
\def\cHH{\enma{\cal H}}
\def\cFF{\enma{\cal F}}
\def\cII{\enma{\cal I}}
\def\cKK{\enma{\cal K}}
\def\cOO{\enma{\cal O}}

\def\LL{\enma{\mathbf{L}}}

\def\boldN{\enma{\mathbf{N}}}
\font\tengoth=eufm10  \font\fivegoth=eufm5
\font\sevengoth=eufm7
\newfam\gothfam  \scriptscriptfont\gothfam=\fivegoth 
\textfont\gothfam=\tengoth \scriptfont\gothfam=\sevengoth

%

\forward{castelnuovo}{Corollary}{4.5}
\forward{small-varieties}{Section}{4}
\forward{restr}{Section}{1}
\forward{sharpness}{Remark}{1.4}
\forward{monomials}{Section}{2}
\forward{Np-bounds}{Section}{3}
\forward{EKS}{Section}{4}
\forward{main1}{Theorem}{1.2}
\forward{main2}{Theorem}{1.1}
\forward{main3}{Theorem}{1.3}
\forward{monomial-Np}{Theorem}{2.1}
\forward{unexpected}{Corollary}{1.7}
\forward{2-regular-gives-small}{Theorem}{1.6}
\forward {monomial schemes lemma}{Proposition}{2.4}
\forward {monomial ideals N}{Corollary}{2.5}


\bigskip
\centerline{\titlefont Restricting linear syzygies: algebra and geometry}
\medskip
\centerline {by}
\smallskip
\centerline {\bf D. Eisenbud, M. Green, K. Hulek, and S. Popescu}

\bigskip

\noindent {\bf Abstract:} 
{\narrower 
\small 
Let $X\subset \PP^r$ be a closed scheme in projective space whose homogeneous 
ideal is generated by quadrics. In this paper we derive geometric consequences
from the presence of a long strand of linear syzygies in the minimal
free resolution of the ideal of $X$.  These consequences are given in
terms of the linear sections of $X$ (the intersections of $X$ with
arbitrary linear subspaces).

More precisely, we say that $X$ (or $I_X$) satisfies $\boldN_{2,p}$ 
if $I_X$ has only linear syzygies for $p$ steps. Thus $X$ is $2$-regular 
in the sense of Castelnuovo-Mumford iff it satisfies $\boldN_{2,p}$ for 
every $p\ge 1$. The simplest of our results says that
if $I_X$ is $2$-regular, then the same is true
for the ideal of any linear section of $X$, 
so long as the intersection has dimension $\leq 1$.
This is not in general true for higher-dimensional
linear sections. We extend this result in
a variety of ways, to projective subschemes satisfying $\boldN_{2,p}$
and to comparisons of resolutions of a subscheme and its linear
sections. We use these results to bound homological invariants 
of some well-known projective varieties.
In Eisenbud-Green-Hulek-Popescu [2004] we use some of the
results of this paper to characterize and classify
all $2$-regular reduced projective schemes.

Extending a result of Fr\"oberg [1990], we give a combinatorial
characterization of the  
monomial ideals  satisfying $\boldN_{2,p}$. Our
results on the $2$-regularity of sections yield a geometric
characterization.

We also apply Green's ``Linear Syzygy Theorem'' [1999] to 
deduce a relation between the resolutions of $I_X$ and 
$I_{X\cup \Gamma}$ for a scheme $\Gamma$, and apply the result 
to bound the number of intersection points of certain pairs of 
varieties such as rational normal scrolls.

}

\bigskip\medskip

\noindent
Let $V$ be a vector space of dimension $r+1$ over an algebraically
closed field $k$ with basis $x_0,\ldots,x_r$.  If
$X\subset\PP^r_k=\PP(V)$ is a nondegenerate closed subscheme we
write $\cII_X$ for the ideal sheaf and $I_X$ for the homogeneous ideal
of $X$ in the homogeneous coordinate ring
$S=\Sym(V)=k[x_0,x_1,\ldots,x_r]$ of $\PP(V)$. Suppose that $I_L$ is an
ideal generated by linear forms, the ideal of a linear space $L$.  In
general there is no strong connection between the minimal free resolution of
$I_X$ and the minimal free resolution of $I_X+I_L$ or of its saturation. 
The goal of this paper is to exhibit some cases where an interesting 
connection of this kind exists.

If $X\subset\PP^r$ is nondegenerate (that is, $I_X$ contains no linear 
form) we say that $X$ {\it satisfies the condition} $\boldN_{2,1}$
if $I_X$ is generated by quadrics. We say that
$X$ {\it satisfies the condition} $\boldN_{2,p}$,
for $p>1$, if in addition the first $p$ steps of the minimal
free resolution 
$$
\cdots \rTo F_t\rTo^{\phi_t} F_{t-1}\rTo^{\phi_{t-1}}\cdots \rTo^{\phi_1} F_0\rTo I_X\rTo 0
$$
of $I_X$ are linear, in the sense that $\phi_t$
is represented by a matrix of linear forms for all 
$1\le t\leq p-1$ or, equivalently, that $\Tor_i^S(I_X,k)$ is a vector space
concentrated in degrees $\leq i+2$ for all $i\leq p-1$.
(Our notation comes from the notation $\boldN_p$ of Green and
Lazarsfeld; but we do not insist that $X$ be projectively
normal, which is their condition $\boldN_0$ and is included in
their condition $\boldN_p$. Also in some of our results we could
replace degree $2$ by an arbitrary degree $d$.) 

In the first section of this paper
we show that if $X\subset\PP^r$ satisfies $\boldN_{2,p}$ then
the same is true of $\Lambda\cap X$ for any
linear subspace $\Lambda$ such that $\dim \Lambda\cap X\leq 1$
and $\dim \Lambda\leq p$. It follows, for example, 
that $\deg (\Lambda\cap X)$ (or even the geometric degree)
is then at most $\dim \Lambda-\dim (\Lambda\cap X)+1$
(see Eisenbud-Green-Hulek-Popescu [2004], Theorem 2.2). 
\ref{main2} may be thought of as a generalization of the
easy direction of Green's conjecture (Green [1984]) proved by Green and
Lazarsfeld [1985], and gives a new proof of this result.
Further, if $\dim \Lambda \leq p-1$, we show that the restriction map 
from the  quadrics in $\PP^r$ containing $X$ to those in $\Lambda$
containing $\Lambda\cap X$ is surjective. As an application
we recover a version of a result of Vermeire [2001] on the linear
system of quadrics through a variety satisfying property
$\boldN_2$. If in addition $X$ is linearly normal, 
$\dim \Lambda\cap X=0$ and $\Lambda\cap X$ spans $\Lambda$,
then the  restriction of minimal free resolutions is surjective for $p-1$ steps. 
We give examples showing that these results are sharp in various
senses.

To describe one of the implications of such results,
we say that a closed subscheme $X\subset\PP^r$ is {\it small\/}
if every zero-dimensional linear section $\Lambda\cap X$ of $X$ 
has $\deg(\Lambda\cap X)\leq 1+\dim \Lambda$. 
If $X$ is nondegenerate, reduced and irreducible, then it follows that
$\deg X=1+\codim(X,\PP^r)$. Such varieties ``of minimal degree''
were classified by Castelnuovo, Del Pezzo, and Bertini (rational
normal curves, scrolls, the Veronese surface, etc). For
reduced subschemes some cases were classified by
by Xamb\`o [1981], and the general case was recently
done by us (see Eisenbud-Green-Hulek-Popescu [2004]) using \ref{main2} from this paper.
In particular, we show there that small algebraic sets are
all $2$-regular. This yields the corollary that if a
reduced subscheme $X\subset\PP^r$ satisfies $\boldN_{2,p}$ for 
$p=\codim(X,\PP^r)$, then $X$ is actually $2$-regular (\ref{unexpected}).

In \ref{monomials} we characterize property $\boldN_{2,p}$
for ideals generated by monomials. In the squarefree case,
an ideal generated by quadratic squarefree monomials comes from a
simplicial complex that is the clique complex of a graph, and
the property $\boldN_{2,p}$ is determined by the length of the shortest
cycle in the graph without a chord (\ref{monomial-Np}; this
result was suggested to us by Serkan Ho\c sten, Ezra Miller,
and Bernd Sturmfels).  A special case is Fr\"oberg's result [1990] 
characterizing $2$-regular square-free monomial ideals. The general monomial case
could be reduced to the squarefree monomial case via polarization, but
we give a direct analysis that yields a more precise result
relating the property $\boldN_{2,p}$ for a monomial ideal to the
corresponding property for the largest squarefree monomial 
ideal it contains (\ref{monomial schemes lemma}).

In a number of interesting cases a kind of strong converse to
\ref{main2} holds: a subscheme $X\subset \PP^r$ satisfies $\boldN_{2,p}$ if and
only if every linear section $\Lambda \cap X$ of dimension zero
satisfies $\deg(\Lambda \cap X)\leq 1+\dim\Lambda$. 
For example Green and Lazarsfeld [1988] prove that this is the case
when $X$ is a smooth linearly normal curve 
of degree $d\geq 3\genus(X)-2$. (See also Eisenbud [2004]
for an exposition and Eisenbud-Popescu-Schreyer-Walter [2002] for
a different perspective.) One of the main results of our paper [2004] shows 
that this is also true for any reduced scheme when $p=\infty$. 
In \ref{monomial ideals N} we  show that this converse holds for any $p\ge 1$ 
when $X$ is defined by a monomial ideal.

In \ref{Np-bounds} we use \refs{main2} and \refn{main1} of \ref{restr}
to prove (conjecturally sharp) upper bounds for the property
$\boldN_p$ for Veronese, Segre-Veronese, Pl\"ucker or Fano embeddings,
as well as for certain embeddings of abelian varieties. In \ref{EKS}
we make use of the Eisenbud-Koh-Stillman conjecture (proved by Green
[1999]) to analyze the restriction of linear syzygies to other simple
subvarieties of the ambient space whose syzygies are known,
e.g. rational normal curves, scrolls, Veronese surfaces, etc. As
applications we give a new proof to Green's syzygetic Castelnuovo
lemma and obtain bounds on the length of a zero-dimensional
intersection of scrolls or Veronese surfaces (for the latter see also
Eisenbud-Hulek-Popescu [2003]).

\medskip
The starting point for this paper was an unpublished note by William Oxbury.
We are grateful to him, as to Aldo Conca, Mark Haiman, J\"urgen Herzog, Jerzy
Weyman, and Sergey Yuzvinsky, for useful discussions
of this material. We are grateful to Serkan Ho\c sten, Ezra Miller,
and Bernd Sturmfels for sparking our interest in the 
monomial question (and suggesting the answer!) 
The program {\it Macaulay2\/} of Mike Stillman and Dan Grayson 
has been extremely useful to us in understanding what
was true about the properties $\boldN_{2,p}$ --- it is
a great advantage only to try to prove true theorems!

Oberwolfach, BIRS, IPAM and MSRI contributed hospitable settings
in which much of this work was done, and we are grateful to
the DFG and the NSF for hospitality and support during the preparation
of this work.

\section{restr} Restricting syzygies to linear subspaces

In this section we show how the condition
$\boldN_{2,p}$ influences low-dimensional linear sections.
 
\theorem{main2} 
Let $X\subset\PP^r$ be a closed subscheme satisfying the property
$\boldN_{2,p}$ with $p\geq 1$, and let $\Lambda\subset\PP^r$ be a linear
subspace of dimension $\leq p$. If $\dim X\cap\Lambda\leq 1$ then
$\cII_{X\cap\Lambda,\Lambda}$ is 2-regular.  In particular, if $X\cap
\Lambda$ is finite, then $\length (X\cap\Lambda)\le \dim\Lambda+1$.

The next two results strengthen \ref{main2} in different ways.
The first is also proven, in a different way, in Eisenbud-Huneke-Ulrich [2004].
A special case of the first follows also from Caviglia [2003].

\theorem{main1} 
Let $X\subset\PP^r$ be a closed subscheme satisfying the property 
$\boldN_{2,p}$ with $p\geq 1$, and let $\Lambda\subset\PP^r$ be a linear
subspace of dimension $\leq p-1$. If $\dim X\cap\Lambda=0$, then the
natural restriction $H^0(\cII_X(2))\to
H^0(\cII_{X\cap\Lambda,\Lambda}(2))$ is surjective.

\theorem{main3} 
Let $X\subset\PP^r$ be a closed subscheme satisfying the property 
$\boldN_{2,p}$ with $p\geq 1$, and let $\Lambda\subset\PP^r$ be a linear
subspace of dimension $\leq p$.  If $X$ is linearly normal,
$X\cap\Lambda$ is zero-dimensional and $X\cap \Lambda$ spans
$\Lambda$, then the natural restriction from the minimal free
resolution of $\cII_{X}$ to the minimal free resolution of
$\cII_{X\cap\Lambda,\Lambda}$ surjects on the first $p-1$ steps.

\remark{sharpness} The conditions in \refs{main2}-\refn{main3}
are often sharp. Here are some examples: 

$1)$ The ideal $I\subset k[x_0,\ldots,x_4]$ of $2\times 2$-minors of
the matrix
$$
\pmatrix{
x_0 & x_1 & 0   & x_2\cr
0   & x_0 & x_1 & x_3}
$$
is saturated and defines a scheme $Y\subset\PP^4$ consisting of a 
$2$-plane with a certain multiplicity $3$ embedded point. 
The scheme $Y$ is a linear section of a $2$-regular variety $X\subset\PP^8$, the
cone over the Segre embedding of $\PP^1 \times\PP^3\subset\PP^7$,
which is $2$-regular and thus satisfies $\boldN_{2,p}$ for every
$p\ge 1$. If $Y$ were at most $1$-dimensional then we would conclude 
from \ref{main2} that $I$ was $2$-regular. However $I$ 
is not even linearly presented. This  shows that the hypothesis 
$\dim(X\cap\Lambda)\leq 1$ in \ref{main2} cannot be weakened.\hfill\break\indent
$2)$ The intersection of $Y$ with the hyperplane $H=\{x_4=0\}$ is $1$-dimensional,
and thus $2$-regular by \ref{main2}. If $Y$ were zero-dimensional we
could conclude from \ref{main1} that the quadrics on $H$
vanishing on $Y\cap H$ were all restrictions of quadrics on $\PP^7$
vanishing on $X$. But the saturation $J$ of $I+(x_4)/(x_4)$ has an 
extra quadratic generator. Thus the hypothesis $\dim(X\cap \Lambda)\leq 0$ in 
\ref{main1} cannot be weakened.\hfill\break\indent
$3)$  The homogeneous ideal 
$$
I=(x_0^2, x_0x_1, x_0x_2-x_1x_4, x_0x_4, x_1x_2-x_1x_4, x_2^2, x_2x_4)\subset k[x_0,\ldots, x_4]
$$
is saturated, satisfies $\boldN_{2,2}$, and defines a scheme $X\subset\PP^4$ 
consisting of two lines concurrent at $p=(0:0:0:1:0)\in\PP^4$ and having an embedded component
at that point. The linear subspace $\Lambda=\{x_3=x_4=0\}$ meets
$X$ in the simple point $q=(0:1:0)$ which does not span $\Lambda$. The
truncation $J$ in degrees $\ge 2$ of the saturation of $I+(x_3,x_4)/(x_3,x_4)$
is thus $2$-regular, but the natural restriction of linear syzygies between
the minimal free resolutions of $I$ and $J$ is not surjective on $\Tor_1$'s.  
Thus the hypothesis $X\cap \Lambda$ spans $\Lambda$ 
in \ref{main3} cannot be weakened.\hfill\break\indent
$4)$ The homogeneous ideal 
$$
I=(x_0^2, x_0x_1-x_2x_4, x_0x_2-x_2x_4, x_0x_3, x_0x_4, x_3x_4, x_4^2)\subset k[x_0,\ldots, x_4]
$$
is the saturated ideal of a $2$-regular scheme $X\subset\PP^4$ consisting of a 
$2$-plane $\Pi$ with two embedded points. Its restriction to the hyperplane 
$\{x_4=0\}$ (which contains the $2$-plane) is a non-saturated ideal 
defining $\Pi$. Its saturation and truncation in degrees $\ge 2$ is a
$2$-regular ideal $J\subset  k[x_0,\ldots, x_3]$, but the restriction map from
the minimal free resolution of $I$ to that of $J$ is not onto. 
This shows that \ref{main3} is sharp.\hfill\break\indent

\remark{d-linear}
Generalizing the condition $\boldN_{2,p}$  we say that a projective subscheme 
$X\subset \PP^r$ {\it satisfies the condition} $\boldN_{d,p}$, for some $d\ge 2$,
if $\Tor_t^S(I_X,k)$ is concentrated in degrees $\leq d+t$
for all $t\leq p-1$. For example, $X$ satisfies $\boldN_{d,1}$ if
$I_X$ is generated in degrees $\leq d$ or, equivalently,
if the truncation $(I_X)_{\geq d}=\oplus_{e\geq d}H^0(\cII_X(e))$
of $I_X$ in degrees $\geq d$ is generated in degree $d$. In
general, it is easy to show that  $X$ satisfies $\boldN_{d,p}$ 
if and only if $X$ satisfies $\boldN_{d,1}$ and
the minimal free resolution $(I_X)_{\geq d}$ is linear
for $p$ steps in the sense above.  The proofs below of \refs{main2} and 
\refn{main1} can be adapted to the case of ideals satisfying 
property $\boldN_{d,p}$, that is generated by forms of any degree $d$ and 
having minimal free resolution with $p$ linear steps (just replace 
in the twists by $\cOO_\Lambda(2)$ with twists of 
$\cOO_\Lambda(d)$).

\bigskip

For the proofs of all three theorems we use the hypercohomology
spectral sequences for the complex obtained by restricting to
$\Lambda$ appropriate twists of the minimal free resolution of the
ideal sheaf $\cII_X$. To fix
notations we recall that if
$$
\cFF^\bullet:\quad\cdots\rTo\cFF^{-m}\rTo\cFF^{1-m}\rTo\cdots\rTo\cFF^{-1}\rTo\cFF^0
$$
is a complex on $\Lambda$, then its hypercohomology $\HH(\cFF^\bullet)$
is computed by two spectral sequences  associated to a Cartan-Eilenberg 
resolution  (double complex) of $\cFF^\bullet$. The filtration by columns 
of the double complex induces a first spectral sequence with $E_2$ terms
$$
\Eprime^{i,j}_2= H^{i}(H^j(\Lambda,\cFF^\bullet))\rImplies \HH^{i+j}(\cFF^\bullet),
$$ 
while the filtration by rows induces a second spectral sequence with
$$
\Edprime^{i,j}_2= H^j(\Lambda,\cHH^{i}(\cFF^\bullet))\rImplies \HH^{i+j}(\cFF^\bullet),
$$ 
where $\cHH^m(\cFF^\bullet)$ denotes the $m$-th cohomology sheaf of
the complex $\cFF^\bullet$.

We start with the case which is technically the simplest to handle:

\bigskip
\noindent{\sl Proof of \ref{main1}:} Let 
$$
\cdots\rTo\cEE^{-n}\rTo\cEE^{-n+1}\rTo\cdots\rTo\cEE^{-1}\rTo\cEE^0\rTo\cII_X\rTo 0
$$
be the sheafification of 
a minimal free resolution of the homogeneous
ideal of $X$. We apply the spectral sequences above to the 
complex 
$\cFF^\bullet:=\cEE^\bullet\otimes\cOO_\Lambda(2)$
obtained by restricting the resolution to $\Lambda$. 

Using the fact that  $X\cap \Lambda$
is zero-dimensional we first show that 
$\HH^{0}(\cFF^\bullet)=H^0(\cII_X\otimes\cOO_\Lambda(2))$.
Since  $\cEE^\bullet$ is a resolution, the sheaves 
$\cHH^i(\cFF^\bullet)$ for $i\le -1$ have support on
the zero-dimensional scheme $X\cap\Lambda$. Hence
$H^j(\Lambda,\cHH^{-i}(\cFF^\bullet))=0$ for all $j\ge 1$ and $i\leq -1$.
Thus the second hypercohomology spectral sequence degenerates at
$\Edprime_2$ and  $\Edprime_2^{i,-i}=0$ for all $i\le 0$. This shows that 
$\HH^{0}(\cFF^\bullet)=\Edprime^{0,0}_\infty = \Edprime^{0,0}_2$.
But $\Edprime^{0,0}_2=H^0(\cII_X\otimes\cOO_\Lambda(2))$
as required since $\cEE^\bullet$
is a resolution of $\cII_X$.

We next use the hypothesis that $X$ satisfies the 
$\boldN_{2,p}$ property to show that the natural restriction map 
from $H^0(\cII_X(2))$ surjects onto $\HH^0(\cFF^\bullet)$, 
which by the result of the previous paragraph
is $H^0(\cII_X\otimes\cOO_\Lambda(2))$.
Consider for this the  other spectral sequence.
By hypothesis $\cFF^i$ is a direct sum of copies of
$\cOO_\Lambda(i)$ for all $1-p\leq i\leq 0$. Since $\dim \Lambda \leq p-1$
$$
\Eprime_1^{i,j}=H^j(\Lambda,\cFF^i)=0
\quad {\rm for}\quad j\ge 1\ {\rm and}\ -\dim \Lambda\le i\le 0.
$$
In particular 
$\HH^0(\cFF^\bullet)=\Eprime^{0,0}_\infty$. Because $\cFF^i\neq 0$ 
only for $i\leq 0$ we see that $\Eprime^{0,0}_1$ surjects via the
natural map onto $\Eprime^{0,0}_\infty$. On the other hand
$\Eprime_1^{0,0}=H^0(\Lambda,\cFF^0)=H^0(\cII_X(2))$ since
$\cEE^\bullet$ is the sheafification of the minimal free resolution
of the homogeneous ideal of $X$.  Combining these maps gives the 
desired surjection.

To complete the proof of the theorem we still need to show that the
natural restriction map
$$H^0(\cII_X\otimes\cOO_\Lambda(2))\rTo H^0(\cII_{X\cap\Lambda,\Lambda}(2))$$
is surjective.
Consider the short exact sequence
$$\leqno{(*)}\qquad\qquad
0\rTo^{}_{}
{\cII_X\cap\cII_\Lambda/ \cII_X\cdot \cII_\Lambda}
\rTo  
\cII_X\otimes\cOO_\Lambda
\rTo 
\cII_{X\cap\Lambda,\Lambda}
\rTo 0.
$$
The kernel $\cKK:=({\cII_X\cap\cII_\Lambda/\cII_X\cdot \cII_\Lambda})$ 
has support on the zero-dimensional scheme $X\cap\Lambda$ so
$H^1(\cKK(2))=0$ and the surjectivity follows.\Box

\medskip

\noindent{\sl Proof of \ref{main2}:}
We use the same two hypercohomology spectral sequences, 
applied this time to each of the complexes
$\cFF^\bullet=\cEE^\bullet\otimes\cOO_\Lambda(2-l)$ for $l\geq 1$.
Using the spectral sequence $\Eprime^{i,j}$ we first prove that 
$\HH^l(\cFF^\bullet)=0.$ To see this note that for
all $i$ with $-(\dim\Lambda-l)\leq i\leq 0$ the term $\cFF^i$
is a direct sum of copies of $\cOO_\Lambda(-l+i)$.
For such $i$ we have
$\Eprime^{i, l-i}_1=H^{l-i}(\Lambda, \cFF^i) = 0$
and the required vanishing follows.

Next we use the second spectral sequence to show that
$$
H^l(\cII_X(2-l)\otimes \cOO_\Lambda) = \Edprime^{0,l}_2 = 0.$$
{}From  $\HH^l(\cFF^\bullet)=0$ it follows that
$\Edprime^{0,l}_\infty = 0$.
The vanishing of the terms 
$$\Edprime_2^{i, -i+l+1}=H^{-i+l+1}(\Lambda, \cHH^i(\cFF^\bullet))$$
is automatic because $\cHH^i(\cFF^\bullet)$
is supported on $X\cap\Lambda$, which has dimension at most
$1< l+1\leq -i+l+1$, for $i\leq 0$. This
implies $\Edprime^{0,l}_2 = \Edprime^{0,l}_\infty = 0$.

Finally, we show that $H^l(\cII_{X\cap\Lambda,\Lambda}(2-l))=0$
by showing that 
$H^l(\cII_X\otimes \cOO_\Lambda(2-l))$ surjects onto it.
For this we use the short exact sequence $(*)$ twisted by $-l$. We
see that it is enough to prove $H^{l+1}(\cKK(2-l))=0$.
This is  automatic because $l\geq 1\geq \dim(X\cap \Lambda)$.\Box

\medskip

\noindent{\sl Proof of \ref{main3}:}
This time we use the spectral sequences on the complex
$\cFF^\bullet=\cEE\otimes \Omega_\Lambda^{m+1}(m+2)$, with
$0\leq m\leq p-1$.
Recall from Green [1984], or Green-Lazarsfeld [1988], or Eisenbud [2004] 
that if $Y\subset \PP^m$ is a scheme with $H^1(\cII_Y(1))=0$, then
for all $m\ge 0$ we have
$$
\Tor_{m}^S(I_Y,k)_{m+2}=H^1(\cII_{Y}\otimes\Omega_{\PP^m}^{m+1}(m+2))
$$
where $S=S_{\PP^m}$ is the homogeneous coordinate ring of $\PP^m$.
Since we have assumed that $X$ is linearly normal we can apply
this to $X\subset \PP^r$. Since 
$X\cap\Lambda$ is $2$-regular by \ref{main2},
and $X\cap \Lambda$ spans $\Lambda$,
we can also apply this with $Y=X\cap \Lambda$ and $\PP^m=\Lambda.$
This gives 
$$
\Tor_{m}^{S_\Lambda}
(I_{X\cap\Lambda,\Lambda},k)_{m+2}=
H^1(\cII_{X\cap\Lambda,\Lambda}\otimes\Omega_{\Lambda}^{m+1}(m+2)).
$$
In the sequence $(*)$ the sheaf $\cKK$ has zero dimensional
support, and we deduce that
$H^1(\cII_{X\cap\Lambda,\Lambda}\otimes\Omega_{\Lambda}^{m+1}(m+2))
=H^1(\cII_{X}\otimes\Omega_\Lambda^{m+1}(m+2))$.

Now consider the spectral sequence $\Edprime$.
We have 
$\Edprime_2^{i,j}=0$ when $i<0$ and $j>0$.
On the other hand we have
$$\Edprime_2^{0,1}=\Tor_{m}^{S_\Lambda}(I_{X\cap\Lambda,\Lambda},k)_{m+2}
$$
by the argument above. For any $q\geq 2$ we have
$$\Edprime_2^{0,q}=H^q(\Lambda, \cII_X\otimes \Omega_\Lambda^{m+1}(m+2)).
$$
These terms are equal to zero because the map
$$\cII_X\otimes \Omega_\Lambda^{m+1}(m+2)
\to
\cOO_{\PP^r}\otimes \Omega_\Lambda^{m+1}(m+2)
$$
has zero-dimensional kernel and cokernel, and 
$H^q(\cOO_{\PP^r}\otimes \Omega_\Lambda^{m+1}(m+2))=0$.
This shows that 
$$
\HH^1(\cFF^\bullet) = 
\Tor_{m}^{S_\Lambda}(I_{X\cap\Lambda,\Lambda},k)_{m+2}.
$$
Next we turn to $\Eprime$. We have 
$\Eprime_1^{i,j}=H^j(\Lambda, \cEE^i\otimes\Omega_\Lambda^{m+1}(m+2))$.
If $0<j<\dim \Lambda$, then Bott's formula gives
$\Eprime_1^{i,j}=0$ unless $j=m+1$ and $i=-m$.
Because $X$ satisfies property $\boldN_{2,p}$ and $m\leq p-1$, we get
$\cEE^{-m} = \Tor_m^S(I_X,k)_{m+2} \otimes \cOO_{\PP^r}(-m-2)$
so
$$
\Eprime_1^{-m,m+1}=H^{m+1}(\Lambda, \cEE^{-m}\otimes\Omega_\Lambda^{m+1}(m+2))
=\Tor_m^S(I_X,k)_{m+2}.
$$ On the other hand if $i\geq -\dim \Lambda +1$ 
then
$
\Eprime_1^{i,\dim \Lambda}
=H^{\dim \Lambda}(\Lambda, \cEE^i\otimes\Omega_\Lambda^{m+1}(m+2))
=0.
$
Thus $\Eprime^{-m,m+1}_1$ surjects onto
$$
\Eprime^{-m,m+1}_\infty =
\HH^1(\cFF^\bullet)=
\Tor_m^{S_\Lambda}(I_{X\cap\Lambda,\Lambda},k)_{m+2}.
$$
This is the natural
map induced by the surjection $I_X\to I_{X\cap\Lambda,\Lambda}$.\Box

Recall from the introduction (see also Eisenbud-Green-Hulek-Popescu [2004]) 
that a closed non-degenerate subscheme
$X\subset\PP^r$ is called ``small'' if every zero-dimensional linear
section of $X$ imposes independent conditions on linear forms, or
equivalently if every zero-dimensional linear section of $X$ is
$2$-regular. 

In this language, we may rephrase \ref{main2} as follows:

\theorem{2-regular-gives-small} Let $X\subset\PP^r$ be a closed subscheme
which satisfies property $\boldN_{2,p}$, for some $p\ge 1$.  If
$\Lambda\subset\PP^r$ is a linear subspace such that
$\codim(X\cap\Lambda, \lspan(X\cap\Lambda))\le p-1$, then $\Lambda\cap
X$ is small.  In particular, $2$-regular projective schemes are
small.\Box

In Eisenbud-Green-Hulek-Popescu [2004] we prove that small and
$2$-regular projective algebraic sets are the same. In particular, we
obtain the following unexpected consequence of
\ref{2-regular-gives-small}:

\corollary{unexpected} If $X\subset \PP^r$
is a reduced subscheme satisfying property $\boldN_{2,p}$
for $p=\codim(X,\PP^r)$, then $X$ is $2$-regular.\Box

For a geometric description and the classification 
of such reduced schemes see Eisenbud-Green-Hulek-Popescu [2004]).

\medskip
As an immediate application of \ref{main1} we get an easy geometric
proof of the known direction of Green's conjecture proved by
Green and Lazarsfeld in [1985].

\corollary{green} Let $C$ be a non-hyperelliptic curve of genus $g\ge 3$,
and let $\nu=\cliff(C)$. Then the canonical embedding of $C$ does not
satisfy property $\boldN_\nu$.

\proof Let $\cOO(L)\in\Pic(C)$ be a degree $d$ line bundle 
that realizes the Clifford index of $C$ with $|L|$ and $|K-L|$
base-point free series. Choose $D_1$ and $D_2$ divisors in $|L|$ and
$|K-L|$, respectively, consisting of distinct points and such that
$D_1\cap D_2=\emptyset$. Since $D_1$ and $D_2$ add up to $K$ their
union spans only a hyperplane and thus the geometric Riemann-Roch
yields that $\Lambda:=\lspan(D_1)\cap\lspan(D_2)$ has dimension
$\nu-1$. On the other hand, $\lspan(D_i)$ meets the canonical curve
only along the points of $D_i$, $i=1,2$, otherwise $C$ would have a
lower Clifford index, and therefore $\Lambda\cap C=\emptyset$.  Thus
from \ref{main1}, if $C$ satisfies property $\boldN_\nu$, we deduce
that $H^0(\cII_C(2))\rTo H^0(\cOO_{\Lambda}(2))$ must be surjective.

Let $D$ be a divisor of degree $\le 2g-2$. From
the cohomology  sequence
$$
\cdots \rTo H^0(\cOO(2K))\rTo H^0(\cOO_D(2K))\rTo H^1(\cOO(2K-D)) \rTo \cdots
$$
we see that $D$ fails to impose independent conditions on
quadrics in $\PP^{g-1}$ if and only if $D\in|K|$. In particular,
$D_1$ and $D_2$ impose independent conditions on quadrics, however
$D_1+D_2$ fails by one to impose independent conditions on quadrics.
This leads now to a contradiction since  the surjectivity
of $H^0(\cII_C(2))\rTo H^0(\cOO_{\Lambda}(2))$, together with
the fact that $D_i$, $i=1,2$, impose independent conditions on quadrics
implies that their sum $D_1+D_2$ would also impose 
independent conditions.\Box

We can also derive a result of Vermeire {\rm [2001]} on rational
mappings of projective space:

\corollary{vermeire} If $X\subset\PP^r$ satisfies
 $\boldN_{2,2}$ and $\Sec(X)\neq\PP^r$, then the linear system
$|H^0(\cII_X(2))|$ on $\PP^r$ is one-to-one outside of $\Sec(X)$.

\proof Let $x_1,x_2\in\PP^r\setminus X$ be a pair of points 
imposing only one condition on the quadrics of $|H^0(\cII_X(2))|$  
and let $\Lambda=\overline{x_1,x_2}$ be the line they span.
By \ref{main1} the restriction map
$H^0(\cII_X(2))\to H^0(\cII_{X\cap\Lambda,\Lambda}(2))$ is surjective and
thus $\Lambda$ must be a secant line of $X\subset\PP^r$.\Box

The following corollary may be regarded as a generalization of \ref{vermeire} 
for the case when property $\boldN_p$, $p\ge 2$ holds:

\corollary{self-associated}
Let $X\subset\PP^r$ be a closed subscheme satisfying the property  
$\boldN_{2,p}$ for some $p\ge 2$, and let $x_1,\ldots,x_p\in\PP^r\setminus X$ 
be points in linearly general position which fail to impose independent 
conditions on the quadrics containing $X$.  Let $\Lambda\cong\PP^{p-1}$
be the linear span of $\{x_1,\ldots,x_p\}$ and assume
that $\Lambda\cap X$ is zero-dimensional and reduced. Then
for some  $2\le q\le p$ there exist subsets $Z_1\subset\{x_1,\ldots,x_p\}$
and $Z_2\subset\Lambda\cap X$, both of cardinality $q$, 
such that $Z_1\cup Z_2$ spans a $\PP^{q-1}$ and fails
(exactly by one) to impose independent conditions on quadrics 
in  $\PP^{q-1}$ (in other words $Z_1\cup Z_2$ is self-associated).

\proof By \ref{main1} the restriction map
$H^0(\cII_X(2))\to H^0(\cII_{X\cap\Lambda,\Lambda}(2))$ is 
surjective, so the hypothesis means that the points $x_1,\ldots,x_p\in\Lambda$
fail to impose independent conditions on the quadrics
in $|H^0(\cII_{X\cap\Lambda,\Lambda}(2))|$. On the other hand, by \ref{main2}
we know that $\deg(X\cap\Lambda)\le p$. The conclusion follows now from
a result of  Dolgachev-Ortland [1988] (Lemma 3, p. 45) and Shokurov [1971] 
which implies that every subscheme of 
$\Gamma:=(\Lambda\cap X)\cup\{x_1,\ldots,x_p\}\subset\Lambda$
of degree $\le 2p$ does impose independent conditions on quadrics in $\Lambda$
if no subset of $2s+2<2p+2$ points of $\Gamma$ is contained in a 
$\PP^s$. (See Eisenbud-Popescu [2000] for the connection with
self-association and the Gorenstein property.)\Box

\section{monomials} Monomial ideals satisfying $\boldN_{2,p}$

In this section we analyze the conditions $\boldN_{2,p}$ for monomial
ideals. We shall see that in the saturated case
(and somewhat more generally) \ref{main2} provides
a criterion to decide which of these conditions are satisfied.

We begin with the case of squarefree monomial ideals.
Using the Stanley-Reisner correspondence, a squarefree
monomial ideal $I\subset S=k[x_0,\ldots,x_r]$ corresponds to a
simplicial complex $\Delta(I)$ with vertices  the variables of the
ring $S$ (see for instance Stanley [1996] for details). We will
denote by $I_\Delta$ the Stanley-Reisner ideal corresponding to a
simplicial complex $\Delta$, and for simplicity we will assume that no
variable $x_i$ is among the minimal generators of $I_\Delta$.

Recall that if $G$ is a graph, then a {\it clique\/} of $G$ is a subset
$T$ of vertices of $G$ such that $G$ contains every edge joining two
vertices of $T$. The {\it clique complex\/} or {\it flag complex\/} of
$G$ is the simplicial complex $\Delta(G)$ whose faces
are precisely the cliques of $G$; the graph $G$ is the
$1$-skeleton of  the clique complex of $G$. 

Clique complexes occur frequently: For any poset $P$ its order complex $\Delta(P)$
(the complex whose faces are the chains of $P$) is the clique complex
$\Delta(G_P)$ of the comparability graph of $P$ (the graph on the
vertices of $P$ whose edges are all pairs of comparable vertices). In
particular the barycentric subdivision of any simplicial complex is 
a clique complex.

It is easy to see that a simplicial complex
$\Delta$ is a clique complex if and only if every minimal
non-face of $\Delta$ consists of $2$ vertices. Thus $I_\Delta$ is generated
by quadratic monomials if and only if $\Delta$ is a clique
complex (of its $1$-skeleton), so to study the properties $\boldN_{2,p}$ we
restrict ourselves to clique complexes. 

The following result was suggested to us by Serkan
Ho\c sten, Ezra Miller, and Bernd Sturmfels. A {\it cycle\/} 
$C$ in $G$ of length $q$ is a sequence of distinct edges of $G$
of the form $(v_1,v_2),(v_2,v_3),\dots,(v_q,v_1)$ joining distinct
vertices $v_1,\ldots,v_q$, for some $q\ge 3$. We say that the 
cycle $C$ has a  {\it chord\/} if some $(v_i,v_j)$ is an edge of $G$,
with $j\not\equiv i+1\ (\hbox{mod }q)$. We say that 
the cycle is {\it minimal\/} if $q>3$ and $C$ has no
chord. Thus the first homology group of $\Delta(G)$
is generated by minimal cycles.

\theorem{monomial-Np} Let $I=I_\Delta$ be
an ideal generated by quadratic squarefree monomials,
and let $G$ be the $1$-skeleton of $\Delta$.
The ideal $I$ satisfies the condition $\boldN_{2,p}$, with $p>1$,
if and only if every minimal cycle in $G$ has length $\geq p+3.$

As a special case we obtain the main theorem of Fr\"oberg [1990].
We say that $G$ is {\it chordal\/} if every cycle of length $>3$
has a chord; in other words if $G$ has no minimal cycles.

\corollary{Froberg}
A square-free monomial ideal
$I_\Delta$ is $2$-regular if and only if $\Delta$ is the
clique complex of a chordal graph.
\Box

The proof of \ref{monomial-Np}
 makes use of Reisner's Theorem (see for example  Hochster [1977], 
or Stanley [1996]):  If $I_\Delta\subset S=k[x_0,\ldots,x_r]$
is a squarefree monomial ideal corresponding to the simplicial complex
$\Delta$, then $\Tor_i^S(I_\Delta, k)$ is a  $\ZZ^{r+1}$-graded vector space 
which is nonzero  only in degrees corresponding to squarefree monomials $m$ and
$$
\Tor_i^S(I_\Delta, k)_m = \widetilde H_{\deg(m)-i-2}(|m|,k),
$$
where $\widetilde H_i(|m|,k)$ denotes the $i$-th reduced homology
of the full subcomplex $|m|$ of $\Delta$ whose vertices 
correspond to the variables dividing $m$. 

\example{elliptic normal curve}
For example, let $d\geq 3$ be an integer. If
$\Delta$ is the simplicial complex with $d+1$ vertices and $d+1$
edges forming a simple cycle, then the reduced homology of
any full proper subcomplex of $\Delta$ is concentrated in degree $0$, while
the reduced homology of the empty set is in degree $-1$ and the
reduced homology of $\Delta$ itself is $k$, concentrated in degree $1$.
We deduce that the minimal free resolution of $S/I_\Delta$
has the form 
$$
0\rTo
S(-d-1)
\rTo 
S(-d+1)^{\beta_{d-2}}
\rTo 
\cdots
\rTo 
S(-2)^{\beta_{1}}
\rTo
S,
$$
so that $I_\Delta$ satisfies $\boldN_{2,d-2}$, but not $\boldN_{2,d-1}$.
Further, the algebraic set $X$ defined by $I_\Delta$ consists of 
$d+1$ lines joined in a cycle in $\PP^d$. The ring $S_X=S/I_\Delta$ is
Cohen-Macaulay, even Gorenstein, and $X$ is a
curve of degree $d+1$ and arithmetic genus $1$ --- a degenerate {\it elliptic normal curve\/}.
If $\Lambda$ is the hyperplane defined by the vanishing of the
sum of the variables (or any hyperplane not containing one
of the components of $X$), then $\Lambda\cap X$ is a set of
$d+1$ points in a $(d-1)$-dimensional plane, and is thus
not $2$-regular.

\medskip
\noindent{\sl Proof of \ref{monomial-Np}.} Let $x_0,\dots,x_r$
be the vertices of $G$, and write $S=k[x_0,\ldots,x_r]$ for
the ambient polynomial ring. Let $X$ be the algebraic set
defined by $I_\Delta$ in $\PP^r$.

First assume that $G$ has a minimal cycle $C$ of length $p+2>3$. 
Let $J$ be the ideal generated by the variables not in the 
support of $C$, and let $\Lambda$ be the projective
linear subspace in $\PP^r$ defined by $J$. The plane section
$\Lambda\cap X\subset\Lambda$ has homogeneous coordinate ring
$S/(I_\Delta+J) = S'/I_C$ where $S'=S/J$. As we showed
in the example above, the ideal $I_C$ is not $2$-regular. 
By \ref{main2}, the ideal $I_{\Delta}$ does not satisfy $\boldN_{2,p}$. 
(Of course the same result may be proven by applying Reisner's Theorem
directly to $\Delta$, by taking $|m|=C$.)

Conversely, suppose that
$I$ does not satisfy the condition $\boldN_{2,p}$,
and take $p>1$ minimal with this property. We must show that
$\Delta$ contains a minimal $(p+2)$-cycle.

By Reisner's Theorem there exists a squarefree
monomial $m$ of minimal degree $\deg(m)\ge p+2$ such that
$\widetilde H_{\deg(m)-p-1}(|m|,k)\neq 0$, while $\widetilde
H_{\deg(m')-i-2}(|m'|,k)=0$ for all $0\le i\le \min(p-1,\deg(m')-3)$ and
all $m'|m$ with $m'\neq m$. If $\deg(m)=p+2$, then $\widetilde H_1(|m|,k)\neq 0$
or equivalently the edge-path group of the simplicial complex $|m|$ is not
trivial. Since $m$ is of minimal degree with the above property, the
simplicial complex $|m|$ must be connected, and again minimality and
the fact that $\Delta$ is a clique complex imply that $|m|$ consists
of a cycle of length $p+2$ in $G$, and this cycle is minimal
 (see also Spanier [1966, Theorem 3, p.140] for a description by generators and relation
of the edge-path group). This is exactly the claim of the theorem.

If however $\deg(m)>p+2$, let $m'|m$ be a squarefree monomial with
$\deg(m')=\deg(m)-1$ and denote by $x$ the extra variable in the
support of $m$.  There is a long exact sequence
$$\ldots\  \widetilde H_i(|{m'}|,k) \rTo \widetilde H_i(|m|,k) \rTo
\widetilde H_{i-1}(\link(x,|m|),k) \rTo \widetilde H_{i-1}(|{m'}|,k)\ \ldots.
$$
which is obtained from the long exact homology sequence of the pair
$(|m|,|m'|)$ and the isomorphisms
$$\widetilde H_i(|m|, |m'|, k) \quad \cong \quad \widetilde H_i(\star(x,|m|),
\link(x,|m|), k) \quad \cong \quad \widetilde H_{i-1}(\link(x,|m|),k)$$
for all $i$. The last isomorphism comes from the
long exact sequence of the second pair which breaks up into isomorphisms
since $\star(x,|m|)$ is contractible. 

Since $\widetilde H_{\deg(m)-p-1}(|m|,k)\neq 0$ while $\widetilde H_{\deg(m)-p-1}(|m'|,k)=0$,
we deduce from the long exact sequence that 
$\widetilde H_{\deg(m)-p-2}(\link(x,|m|),k)\neq 0$, with $\deg(m)-p-2\ge 1$. 
On the other hand the simplicial complex $\link(x,|m|)$ is a full (strict) subcomplex
of $|m|$ and thus of $\Delta$. Indeed if $x_{i_1},\ldots, x_{i_s}\in \link(x,|m|)$ are
vertices such that $\{x_{i_1},\ldots, x_{i_s}\}\in |m|\subseteq \Delta$, then obviously  
$\{x_{i_a}, x_{i_b}\}\in\Delta$ for all $a\neq b$, and also $\{x, x_{i_a}\}\in\Delta$ 
by the definition of the link. Since $\Delta$ is a clique complex it 
follows that $\{x, x_{i_1},\ldots, x_{i_s}\}$ must also be a face of 
$\Delta$ with support in $|m|$. But this means that we  have found a 
full subcomplex $|m''|=\link(x,|m|)$ of $\Delta$, with $\deg(m'')<\deg(m)$, 
such that $\widetilde H_j(|m''|, k)\neq 0$ for some $j\ge 1$, which 
contradicts the fact that $I_\Delta$ satisfies property $\boldN_{2,p-1}$.
This concludes the proof of the theorem.\Box

As Fr\"oberg remarks, the case of a general ideal $I\subset S=k[x_0,\dots,x_r]$ 
generated by quadratic monomials may be reduced, by the process of polarization,
to the squarefree case. However, we can give a more explicit result. 
We may harmlessly assume that $I$ contains no linear forms, and we
may write $I$ uniquely in the form
$I=I_\Delta+I_s$ for some simplicial complex
$\Delta$ with vertices $x_0,\dots,x_r$
and where the ideal $I_s$ is generated by $\{x_i^2\mid x_i^2\in I\}$. We will refer
to the vertices $x$ of $\Delta$ such that $x^2\in I$
as the {\it square vertices for $I$\/}. 

\proposition{monomial schemes lemma}
Let $I=I_\Delta+I_s$ be an ideal generated by
quadratic monomials,  decomposed as above.
\item{$a)$} The ideal $I$ satisfies $\boldN_{2,2}$ if and
only if $I_\Delta$ satisfies $\boldN_{2,2}$ and for any square vertex 
$x$ for $I$, $\link(x,\Delta)$ is a simplex not containing any 
square vertex for $I$.
\item{$b)$} If $I$ satisfies $\boldN_{2,2}$, then $I$ satisfies
$\boldN_{2,p}$ for some $p\ge 3$ if and only if $I_\Delta$ 
satisfies $\boldN_{2,p}$.

\proof If $I_s=(0)$ the result is obvious. Otherwise,
let $x$ be a square vertex for $I$, and let 
$I'=I_\Delta+I_s'$, where $I_s'\subset I_s$ is the ideal generated by
the squares of all square vertices for $I$ other than $x$. The 
exact sequence
$$
0\rTo ((I':x^2)/I')(-2)\rTo S/I'(-2) \rTo^{x^2}S/I'\rTo S/I\rTo 0
$$
and the observation that $(I':x^2)=(I':x)$ yields a short exact
sequence
$$
0\rTo(S/(I':x))(-2) \rTo^{x^2}S/I'\rTo S/I\rTo 0.
$$
{}From the long exact sequence in $\Tor$'s,
 we see that $I$ satisfies property $\boldN_{2,2}$ if and only
if $I'$ satisfies $\boldN_{2,2}$ and $(I':x)$ is generated by
linear forms. On the other hand, we have $(I':x)= I_{\link(x,\Delta)} + I'_s$.
This is generated by linear forms if and only if
$\link(x,\Delta)$ is a simplex not containing any of the 
square vertices that appear in $I'_s$, as required. This proves part $a)$.

When $(I':x)$ is generated by linear forms, each
$\Tor_i^S(S/(I':x), k)$ is concentrated in degree $i$.
In this circumstance the long exact sequence 
in $\Tor$'s coming from the short exact sequence
above shows that $I$ satisfies $\boldN_{2,p}$ for some $p\ge 3$ if and only
if $I'$ satisfies $\boldN_{2,p}$, and we are done by induction. 
\Box

\corollary{monomial ideals N} If $I=I_X$ is the ideal
of a closed subscheme $X\subset\PP^r$, and $I$ is  
generated by quadratic monomials, then $I$ satisfies
$\boldN_{2,p}$ if and only if for all planes $\Lambda$ of
dimension $\leq p$ having zero-dimensional intersection
with $X$ the scheme $\Lambda\cap X$ is $2$-regular.

We first need a characterization of saturated ideals:

\lemma{saturated monomials} Let $I=I_\Delta+I_s$ be an
ideal generated by quadratic monomials, decomposed as above,
with $I_\Delta$ a squarefree quadratic monomial ideal
and $I_s$ the ideal generated by the squares of 
the square vertices for $I$. Then 
$I$ is saturated if and only if every maximal face of
$\Delta$ contains at least one non-square vertex for $I$.

\proof If the ideal generated by all the vertices is 
associated, it must annihilate a squarefree monomial,
and this must be the product of all vertices in a facet of $\Delta$. 
Such a product is annihilated by the maximal ideal if and 
only if every vertex in that facet is
a square vertex for $I$.\Box

\noindent{\sl Proof of \ref{monomial ideals N}.\/}
If a linear subspace $\Lambda$ of dimension $\leq p$ meets
$X$ in a zero-dimensional scheme $X\cap\Lambda$ that is not $2$-regular, 
then \ref{main2} shows that $I$ does not satisfy $\boldN_{2,p}$.

Conversely, suppose that $I$ does not satisfy $\boldN_{2,p}$, with
$p\ge 2$ minimal, and decompose $I=I_\Delta+I_s$ as above.  If $p>2$
then from \ref{monomial schemes lemma} $b)$ we see that $I_\Delta$
does not satisfy $\boldN_{2,p}$, and thus the $1$-skeleton of $\Delta$
has a minimal cycle $C$ of length $p+2$. If $x$ is a vertex of such a
cycle then $\link(x,\Delta)$ is not a simplex, and it follows that $x$
is not a square vertex for $I$.  If $\Lambda'$ is the linear subspace
spanned by all the vertices in the cycle $C$, then $X\cap
\Lambda'\subset\Lambda'$ is a degenerate ``elliptic normal curve'' as
in \ref{elliptic normal curve}.  As remarked in that example, any
sufficiently general plane $\Lambda\subset \Lambda'$ of codimension
$1$ in $\Lambda'$ is a $p$-plane that meets $X$ in a zero-dimensional
scheme that is not $2$-regular.

Finally, suppose that $I$ does not satisfy $\boldN_{2,2}$.
We use the characterization in part $a)$ of \ref{monomial schemes lemma}.
If $I_\Delta$ does not satisfy $\boldN_{2,2}$ then we proceed as before.
Otherwise there is a square vertex $x$ for $I$ such that either the
link of $x$ in $\Delta$ is not a simplex, or the link of $x$
in $\Delta$ is a simplex containing another square vertex
for $I$. 

In the first case we can choose vertices $y,z$ in $\link(x,\Delta)$ 
such that $yz\in I_\Delta\subseteq I$. Factoring out all the variables
except $x,y,z$ we get from $I$ a monomial ideal $\overline I$ with
$$
(x^2, yz)\subset \overline I\subset (x^2, y^2, yz, z^2)\subset k[x,y,z].
$$
Any such ideal is saturated, so it defines the zero-dimensional 
scheme $X\cap \Lambda\subset\Lambda$, where $\Lambda$ is the $2$-plane 
spanned by the vertices $x,y,z$. Further, $\overline I$ is not $2$-regular.

In the second case, let $y$ be one of the square vertices for
$I$ such that $y\in\link(x,\Delta)$. Since the link of $x$ is 
a simplex, the star of $x$ (that is, $x$ together with the link)
is a maximal face of $\Delta$. Since $I$ is saturated we may choose
a vertex $z\in\link(x,\Delta)$ that is not a square vertex 
for $I$. Factoring out all the variables except $x$, $y$
and $z$ we get from $I$ the saturated ideal $\overline I=(x^2, y^2)\subset k[x,y,z]$,
an ideal that is not $2$-regular, and that is the ideal of
a zero-dimensional intersection of $X$ with a $2$-plane. This
concludes the proof of the corollary.
\Box

\corollary{independence-of-char}
The condition that a monomial ideal in $S=k[x_0,\dots,x_r]$ satisfies
property $\boldN_{2,p}$ for some $p\ge 1$, and in particular $2$-regularity, 
is independent of the field $k$ (not necessarily
algebraically closed).\Box

\remark{initial-ideal-remark} By a result of Bayer and
Stillman [1987] (or see Eisenbud [1995] Theorem 15.20)
a subscheme $X\subset \PP^r$ over a field of characteristic
zero is $2$-regular if and only if it has a Borel-fixed (generic) 
initial ideal generated by quadratic monomials. Any scheme 
defined by a monomial ideal is, moreover, the degeneration by a flat family of 
linear sections, of a reduced union $Y$ of planes defined by the 
monomials of a ``polarization'' (see for example Eisenbud
[2004]). Thus each $2$-regular
scheme $X$ over a field of characteristic zero, reduced or not, 
is associated canonically with an absolutely reduced
scheme $Y$, a union of coordinate planes, that is also $2$-regular. 

\example{} If $X$ is the union of two disjoint lines in 
$\PP^3$ with ideal $(a,b)\cap (c,d)\subset k[a,b,c,d]$
then $X$ has generic initial ideal 
$$
(x^2, xy, y^2,xz)\subset k[x,y,z,w]=k[a,b,c,d],
$$
and this ideal has polarization 
$$
(x_1x_2, x_1y_1, y_1y_2, x_1z)\subset k[x_1,x_2,y_1,y_2,z, w].
$$
In this case the polarization scheme $Y\subset\PP^5$ is the 
cone over the reduced union of two planes in $\PP^3$ 
meeting a line ``sticking out into'' $\PP^4$ in a point 
on the intersection of the two planes.
Thus even if the original scheme is a union of planes,
the resulting polarization of its generic initial 
ideal may be quite different. In this case the 
general hyperplane section of $Y$ is the cone over the
scheme consisting of $2$ lines in $\PP^2$ with
a spatial embedded point of multiplicity one at their intersection --
the limit of the original scheme $X$ in a family where the
two lines become coplanar.

\section{Np-bounds} Upper bounds for property $\boldN_p$

{}From an alternative perspective the results in \ref{restr} provide
geometric explanations for the failure of property $\boldN_p$ and thus
allow one to test optimality of the results of Green, Ein-Lazarsfeld,
and many others, mentioned in the introduction.

Perhaps the simplest example (handled by different methods in
Ottaviani-Paoletti [2001]) is the necessity of the conditions in the
following:

\conjecture{veronese-conjecture} 
Property $\boldN_p$ holds for the 
$d$-uple embedding of $\PP^n$ if and only if
either 
\item{$-$}$n=1$ and $d,p\in\NN$, or 
\item{$-$}$n=2$, $d=2$ and $p\in\NN$, or
\item{$-$}$n\ge 3$, $d=2$, and $p\le 5$, or 
\item{$-$}$n\ge 2$, $d\ge 3$ and $p\le 3d-3$.

Jozefiak-Pragacz-Weyman [1981] show that the $2$-uple embedding of
$\PP^n$, $n\ge 3$, satisfies property $\boldN_{5}$. In the case of the
$d$-uple embedding of $\PP^2$ its minimal free resolution restricts
to the minimal free resolution of a hyperplane section (a plane curve),
and so Green [1984] implies that for $d\ge 3$ the $d$-uple
embedding of $\PP^2$ satisfies property $\boldN_{3d-3}$. See also
Rubei [2003] for a proof of the fact that $3$-uple embedding of $\PP^n$
satisfies property $\boldN_{4}$ for all $n$. In all other
cases the sufficiency of the conditions in \ref{veronese-conjecture}
is wide open. On the other hand \ref{main2} yields easily the
necessity of those conditions, namely

\proposition{veronese} Let $n\ge 2$ and $d\ge 2$ be integers.
\item{$(a)$} If $n\ge 2$ and $d\ge 3$, then the $d$-uple embedding 
of $\PP^n$ fails property $\boldN_{3d-2}$.
\item{$(b)$} If $n\ge 3$, then the $2$-uple embedding of 
$\PP^n$ fails property $\boldN_{6}$.

\proof Observe first that for all $m<n$ the $d$-uple embedding of 
$\PP^m$ is a linear section of the $d$-uple  embedding of 
$\PP^n$ and thus by \ref{main2} for the failure of property $\boldN_p$ 
it is enough to produce a $(p+2)$-secant $p$-plane to the $d$-uple 
embedding of $\PP^m$ for some $m<n$. 

To prove $a)$ we may assume that $n=2$. Since $d\ge 3$, a complete intersection 
$(3,d)$ in $\PP^2$ is cut out  by forms of degree $d$ but fails to impose 
independent conditions on such forms. In other words, the linear span of the 
$d$-uple embedding of such a complete intersection is a $3d$-secant linear space 
of dimension $3d-2$ to the $d$-uple embedding of $\PP^2$, which therefore 
must fail property $\boldN_{3d-2}$ by \ref{main2}.

Similarly, a complete intersection $Z\subset\PP^3$ of three quadrics
fails to impose independent conditions on quadrics, and thus maps via
the $2$-uple  embedding of $\PP^3$ to a collection of $8$
points spanning only a six dimensional linear subspace $\Lambda$ of
$\PP^9$ meeting the $2$-uple embedding of $\PP^3$
only along the points of $Z$. By \ref{main2}, the $2$-uple Veronese
embedding of $\PP^3$ fails property $\boldN_{6}$.\Box

The failure of property $\boldN_{3d-2}$ for the $d$-uple embedding of
$\PP^2$ can be accounted for also by the existence of a relatively
long strand of linear syzygies in the minimal free resolution of
$\omega_{\PP^2}(d)$.  Namely, with notations as in
Eisenbud-Popescu-Schreyer-Walter [2002], we have the following

\proposition{veronese-BGG} Let $W=H^0(\omega_{\PP^2}(d))$ and 
set $w=\dim(W)$, let $U=H^0(\omega_{\PP^2}^{-1})$, let
$V=H^0(\cOO_{\PP^2}(d))$ and $S=\Sym(V)$.  If $d\ge 3$, the natural
multiplication pairing $\mu:W\otimes U\to V$ makes
$Q=\oplus_l(\wedge^{l+1}(W^*)\otimes\Sym^l(U^*))$ into a graded
$E=\wedge^*(V^*)$-module such that the maximal irredundant quotient of the
linear complex
$$\LL(Q^*):\;
0\to \wedge^w W\otimes D_{w-1}(U)\otimes S(-w+1)\to\cdots
\to\wedge^2 W\otimes U\otimes S(-1)\to W\otimes S
$$
is a linear complex of the same length which injects as a degreewise
direct summand into the minimal free resolution of the $S$-module
$\oplus_{m\ge 0}H^0(\omega_{\PP^2}(d(m+1)))$. In particular the
$d$-uple embedding of $\PP^2$ fails property $\boldN_{3d-2}$.

\proof The above multiplication pairing $\mu$ is obviously geometrically
$1$-generic so the first part of the claim is
a direct application of Proposition 2.10 in 
Eisenbud-Popescu-Schreyer-Walter [2002] with
$L=\cOO_{\PP^2}(d-3)$, $L'=\cOO_{\PP^2}(3)$, and $L''=\cOO_{\PP^2}(d)$,
and with $W$, $U$, and $V$ as in the statement of the
proposition. 

For the second claim observe first that the homogeneous ideal $I_d$ of
the $d$-uple embedding of $\PP^2$ is generated by quadrics but is only
$3$-regular so that its minimal free resolution has two strands
(linear and quadratic).  On the other hand the dual of the maximal
irredundant quotient of the linear complex $\LL(Q^*)$ has length
${{d-1}\choose 2}$ and is a degreewise direct summand of the second
strand into the minimal free resolution of $I_d$. Since the whole
resolution of $I_d$ has length ${{d+2}\choose 2}-3$ it follows that
the $d$-uple embedding of $\PP^2$ fails property $\boldN_{3d-2}$.\Box

The argument used in the proof of \ref{veronese} $a)$ provides
upper bounds for property $\boldN_p$ for other Fano-type
varieties and embeddings. For instance for  
embeddings of ruled or Del Pezzo surfaces we obtain the following bounds 
(where \ref{veronese} $a)$ is the case when $S=\PP^2$): 

\proposition{surface-fano} Let $S$ be a smooth surface and 
$L$ be a very ample line bundle on $S$. If $|-K_S|\ne\emptyset$,
and $\cOO(K_S)\otimes L$ is globally generated, then
the image of $S$ via the linear system $|L|$ fails
property $\boldN_{-K_S\cdot L-2}$.

\proof  Let $D\in |-K_S|$, let $C\in|L|$ be a general
curve and denote by $Z=D\cap C$ their intersection. The Koszul complex
on the sections defining $D$ and $C$ expands to the 
following commutative diagram
$$
\diagram[midshaft,small]
  &     &0          &     &0      &     &0      &     &\cr
  &     &\uTo       &     &\uTo   &     &\uTo   &     &\cr
0 &\rTo &\cOO_C(-D) &\rTo &\cOO_C &\rTo &\cOO_Z &\rTo &0\cr 
  &     &\uTo       &     &\uTo   &     &\uTo   &     &\cr
0 &\rTo &\cOO_S(-D) &\rTo &\cOO_S &\rTo &\cOO_D &\rTo &0\cr
  &     &\uTo       &     &\uTo   &     &\uTo   &     &\cr
0 &\rTo&\cOO_S(-D-C)&\rTo &\cOO_S(-C)&\rTo &\cOO_D(-C)&\rTo &0\cr
  &     &\uTo       &     &\uTo   &     &\uTo   &     &\cr
  &     &0          &     &0      &     &0      &     &\cr
\enddiagram
$$
which we twist by $L$ and take cohomology.  
{}From the long exact sequence of the middle row, since $H^1(\cOO(-D)\otimes L)=
H^1(\cOO(K_S)\otimes L)=0$ by Kodaira vanishing (in characteristic 0) or by
Shepherd-Barron [1991] and Terakawa [1999, Theorem 1.6] (in positive characteristic), 
we deduce that the natural restriction map $H^0(L)\to H^0(L_{\mid D})$ is
surjective. On the other hand in the long exact sequence of the last
column 
$$
\cdots \rTo H^0(L_{\mid D})\rTo H^0(L_{\mid Z})\rTo H^1(\cOO_D(-C)\otimes L)
\rTo H^1(L_{\mid D})\rTo\cdots
$$
we have $h^1(\cOO_D(-C)\otimes L)=h^0(\cOO_D)\ge 1$, while 
$h^1(L_{\mid D})= h^0(L^{-1}_{\mid D})=0$ since $\cOO(K_D)=\cOO_D$
and $L$ is ample. Putting everything together it follows that  
the subscheme $Z$ fails to impose independent conditions
on the sections of $L$. Moreover since $\cOO(K_S)\otimes L$
is globally generated we deduce that $\cII_Z\otimes L$ is 
also globally generated. But $\length(Z)= -K_S\cdot L$ so
the claim of the proposition follows now directly from
\ref{main2}.\Box

\remark{} $1)$ By adjunction (see Sommese [1979] or Sommese-Van de Ven [1987]),
in \ref{surface-fano}, the line bundle $\cOO(K_S)\otimes L$ is globally
generated if and only if $(S,L)$ is not one of the following pairs:
$(\PP^2,\cOO_{\PP^2}(1))$, $(\PP^2,\cOO_{\PP^2}(2))$, or
$(\PP(E),\cOO_{\PP(E)}(1))$ with $E$ a rank $2$ vector bundle on a
curve.\hfill\break\indent
$2)$ A similar argument as in \ref{surface-fano} shows that if 
$X$ is a smooth projective surface and $L$ is a very ample divisor on it, 
then the embedding of $X$ via the linear system $|K_X+(p+3)L|$ fails to 
satisfy property $\boldN_{3pL^2-2}$, for $p\ge 3$ (or fails to satisfy property
$\boldN_{(2p+2)L^2-2}$ for $p\ge 2$, if $(X,\cOO(L))\neq(\PP^2,\cOO_{\PP^2}(1))$).

\proposition{segre} Let $X$ denote the image of the
Segre-Veronese embedding
$$
\PP^{n_1}\times\PP^{n_2}\times\cdots\times\PP^{n_m}
\rInto^{\quad(d_1,d_2,\ldots,d_m)\quad}\PP^{\prod_{i=1}^m{{n_i+d_i}\choose d_i}-1}
$$
\item{$(a)$} If $m\ge 3$ and $d_i=1$ for at least three values of $1\le i\le m$, 
then  $X$ fails property $\boldN_4$.
\item{$(b)$} If $m\ge 3$ and $d_i=1$ for exactly two values of  $1\le i\le m$, then
$X$ fails property $\boldN_{2\min_{\{i\mid d_i>1\}}d_i+2}$. 
\item{$(c)$} If $m\ge 3$ and $d_i=1$ for at most one value of  $1\le i\le m$,
or if $m\ge 2$ and $d_i>1$ for all $1\le i\le m$, then $X$ fails property 
$\boldN_{2\min_{\{i\ne j\mid d_i,d_j>1\}}(d_i+d_j)-2}$.

\proof We may argue as in the proof of \ref{veronese}
and exhibit for suitable $p$ a $p$-dimensional linear
subspace which is $(p+2)$-secant to the Segre-Veronese 
embedding of a product of $r<m$ factors. Failure
of property $\boldN_p$ then follows from \ref{main2}.

To prove $a)$ we may assume that $m=3$. The linear span
of the Segre-Veronese embedding of a complete intersection of 
type ${(1,1,1)}^3$ is a $6$-secant $\PP^4$, thus $X$ fails
property $\boldN_4$ in this case.

Case $b)$ is similar to $a)$: we may assume that $m=3$
and consider the linear span of the Segre-Veronese embedding of
a complete intersection of one hypersurface of
multidegree $(1,1,2)$, and two hypersurfaces of
multidegree $(1,1,d)$ with $d=\min_{\{i\mid d_i>1\}}d_i$.

Finally in case $c)$ we may assume that $m=2$ and that both degrees
are $\ge 2$, in which case the claim follows from \ref{surface-fano}
for $S=\PP^1\times\PP^1$. 
\Box

Not much is known concerning the converse of \ref{segre}, except for
computational evidence via Macaulay2 for small values of $m$, $n_i$
and $d_i$. The following remark collects all positive related results
we are aware of. 

\remark{many-ones} $1)$ If $d_1,d_2\ge 2$, then the embedding of
$\PP^1\times\PP^1$ via the linear system 
$|\cOO_{\PP^1\times\PP^1}(d_1,d_2)|$  satisfies
$\boldN_{2d_1+2d_2-3}$ (see Gallego-Purnaprajna [2001]),
but fails to satisfy $\boldN_{2d_1+2d_2-2}$ by \ref{segre} or
\ref{surface-fano} above.\hfill\break\indent
$2)$ Lascoux [1978] and Pragacz-Weyman [1985] describe the minimal free
resolution of the  Segre embedding of $\PP^{n_1}\times\PP^{n_2}$. In
particular they show that it satisfies property $\boldN_p$ if and only 
if $p\le 3$.\hfill\break\indent
$3)$ Using simplicial methods Rubei [2002, 2004] shows that the Segre 
embedding of $\PP^{n_1}\times\PP^{n_2}\times\cdots\times\PP^{n_m}$ 
(at least three factors) satisfies property $\boldN_p$ if and only 
if $p\le 3$. Also Corollary 8 in Rubei [2002] proves part $b)$ 
in \ref{segre} via a different method.\hfill\break\indent
$4)$ The resolution of the Segre embedding of $\PP^1\times\PP^1\times\PP^1$
as well as a number of other special cases where the resolution is
self-dual are investigated in Barcanescu-Manolache [1981].

\proposition{pluecker} The Pl\"ucker embedding of the Grassmannian
$\Gr(k,n)\subset\PP^{{n\choose k}-1}$, where $2\le k\le n-2$ and $n\ge 5$, 
fails property $\boldN_3$.

\proof It is enough to observe that for all $2\le k\le n-2$ and $n\ge 5$,
the Pl\"ucker embedding of the Grassmannian
$\Gr(k,n)\subset\PP^{{n\choose k}-1}$ has as linear section the
Pl\"ucker embedding of $\Gr(2,5)\subset\PP^{9}$. On the other hand a
general codimension three linear section of $\Gr(2,5)\subset\PP^{9}$
is a collection of $5$ points spanning only a $\PP^3$, and the
conclusion follows now again from \ref{main2}.
\Box

\remark{} $1)$ Jerzy Weyman informed us that property $\boldN_2$ 
always holds for the  Pl\"ucker embedding of any Grassmannian.
\hfill\break\indent
$2)$ Manivel [1996] proved that if $X=G/P$, where $G=\SL(V)$, 
$V$ is a complex vector space and $P$ a parabolic subgroup,
and $L$ is a very ample line bundle on $X$, then the embedding 
defined by the complete linear system $|L^p|$ satisfies property 
$\boldN_p$ for all $p\ge 1$.

\medskip

Recall that a complete linear system $|L|$ on a projective variety
$X$ is said to be $k$-very ample if for any zero dimensional
subscheme $Z\subset X$ of length $k+1$ the restriction map
$$
H^0(L)\rTo H^0(L_{\mid Z})
$$
is surjective. In particular $0$-very ample is 
``base point free'' and $1$-very ample is ``very ample''.

Pareschi [2000] and Pareschi-Popa [2003] proved that if $X$ is an
abelian variety and $L_1,\ldots,L_{p+3}$ are ample line bundles on $X$
then the embedding of $X\subset\PP^N$ by the linear system
$|L_1\otimes\ldots\otimes L_{p+3}|$ satisfies property $\boldN_p$. By
\ref{main2}  and the classification of small algebraic sets in
Eisenbud-Green-Hulek-Popescu [2004] every positive dimensional 
reduced irreducible component of a linear
section $\Lambda\cap X$, where $\Lambda$ a linear subspace of dimension
$\le p$ of $\PP^N$, is a variety of minimal degree in its linear span
and hence rational. But abelian varieties do not contain rational
positive dimensional subvarieties. Thus by
\ref{main2} it follows that $L_1\otimes\ldots\otimes L_{p+3}$
is $(p+1)$-very ample, which is a special case of
Theorem 1 in Bauer-Szemberg [1997]). 

 Observe also that if $X=\prod_{i=1}^{\dim(X)} E_i$ is a product of
elliptic curves, each with origin $o_{E_i}$, and $L:=\prod_i
p_i^*(\cOO_{E_i}(o_{E_i}))$ is the canonical principal polarization on
$X$, then $L^{p+3}$ fails to satisfy property $\boldN_{p+1}$.  This is
a consequence of \ref{main2} and Abel's theorem since one may choose
$(p+3)$ points on $E_i$ such that any divisor in the linear system
$|(p+3)o_{E_i}|$ containing $(p+2)$ of those points contains also the
remaining point.

Gross-Popescu [1998] conjectured that the general $(1,d)$-polarized
abelian surface, for $d\ge 10$, satisfies property $\boldN_{\lbrack
{d\over 2}\rbrack-4}$. As above, by \ref{main2}, this would imply that
a $(1,d)$-polarization on a general abelian surface is $k$-very ample
if $d\ge 2k+3$ and $d\ge 10$ (compare again with Bauer-Szemberg
[1997], Theorem 1).

\section{EKS} Secants and syzygy varieties

In this section we analyze the restriction of linear syzygies to
non-linear varieties with known syzygies such as rational normal
curves, rational scrolls and Veronese surfaces.
 
\theorem{secants} Let $X, \Gamma\subset\PP^r$ be subschemes 
such that $X$ is non-degenerate and  $\Gamma$ is reduced with
every irreducible component spanning all of $\PP^r$.
If the natural restriction map
$$
\Tor^S_p(I_{X\cup\Gamma},k)_{p+2}\rTo \Tor^S_p(I_{X},k)_{p+2}
$$
is not surjective, then 
$$
\dim {{H^0(\cII_{X\cap\Gamma}(2))}\over{H^0(\cII_\Gamma(2))}} > p.
$$

\proof The main idea will be to use the Eisenbud-Koh-Stillman
Conjecture {\bf EKS} (proved by M. Green [1999]) which says that if
$M=\oplus_{i\ge 0}M_i$ is a finitely generated graded module over the
polynomial ring $S=\Sym(V)$ such that
\item{$a)$} $\ker(\wedge^pV\otimes M_0\rTo\wedge^{p-1}V\otimes M_1)\ne 0$, for
some $p>0$, and moreover
\item{$b)$} $\dim M_0\le p$,

\noindent
then there exist a $p$ dimensional family of rank one relations
(i.e. decomposable tensors) in the kernel of the multiplication map
$V\otimes M_0\rTo M_1$.

We will not need the full strength of {\bf EKS}, but just the
existence of such rank one relations under the above hypothesis, 
and we will apply {\bf EKS} to
$$
M=\bigoplus_{i\ge 0}{{H^0(\cII_{X\cap\Gamma}(i+2))}\over{H^0(\cII_\Gamma(i+2))}}
$$
regarded as a finitely generated module over $S=\Sym(V)$, the polynomial ring
of the ambient $\PP^r$.

There are no rank one relations in the kernel of the multiplication
morphism $V\otimes M_0\rTo M_1$. Such a rank one relation would amount
to the existence of a quadric defined by $Q\in
H^0(\cII_{X\cap\Gamma}(2))$ not vanishing on $\Gamma$ and a
hyperplane defined by $H\in H^0(\cOO_{\PP^r}(1))$ such that $QH\in
H^0(\cII_{\Gamma}(3))$, which is impossible since each irreducible
component of $\Gamma$ is assumed to be nondegenerate.

We will relate condition $a)$ in {\bf EKS} for the module $M$ to
the analogous one for the module
$$
P=\bigoplus_{i\ge 0}{{H^0(\cII_{X}(i+2))}\over{H^0(\cII_{X\cup\Gamma}(i+2))}}.
$$
Expressing as usual the Tor's via Koszul cohomology our hypothesis
that 
$$
\Tor^S_p(I_{X\cup\Gamma},k)_{p+2}\rTo \Tor^S_p(I_{X},k)_{p+2}
$$
is not surjective translates into the existence of an element
$$\alpha\in\ker(\wedge^p V\otimes H^0(\cII_{X}(2))\rTo \wedge^{p-1} V\otimes 
H^0(\cII_{X}(3)))$$ 
which is not in the image 
of the natural inclusion morphism
$$\wedge^p V\otimes H^0(\cII_{X\cup\Gamma}(2))\rTo \wedge^{p} V\otimes 
H^0(\cII_{X}(2)).$$
Taking global sections in the first row of the exact diagram of ideal sheaves
$$\diagram[small,midshaft]
0&\rTo&\cII_{X\cup\Gamma}(2)&\rTo&\cII_{X}(2)&\rTo&\cII_{X\cap\Gamma,\Gamma}(2)\\
&&\dInto&&\dInto&&\dTo^{\cong}\\
0&\rTo&\cII_{\Gamma}(2)&\rTo&\cII_{X\cap\Gamma}(2)&\rTo&\cII_{X\cap\Gamma,\Gamma}(2)&\rTo&0\\
\enddiagram$$
we see that $\alpha$ induces a non-trivial element $\bar\alpha$ in
$$\bar\alpha\in\ker(\wedge^p V\otimes P_0\rTo \wedge^{p-1} V\otimes P_1).$$ 
On the other hand, twisting and taking global sections in the above diagram
yields the inclusion $P\subseteq M$. In particular, we may view
$\bar\alpha$ as an element of 
$\ker(\wedge^pV\otimes M_0\rTo\wedge^{p-1}V\otimes M_1)$, 
which is thus non-zero. By {\bf EKS}, since there are no rank one
relations in the kernel of $V\otimes M_0\rTo M_1$, we deduce that
$\dim M_0 > p$ which finishes the proof of the theorem.\Box

In the case where $\Gamma$ is a smooth curve \ref{secants} 
has more geometric content and so we restate a special case of 
it explicitly:

\corollary{secants-curve} With notation
as in \ref{secants}, if $\Gamma$ is an irreducible 
nondegenerate curve such that
$$ 
\Tor^S_p(I_{X\cup\Gamma},k)_{p+2}\rTo \Tor^S_p(I_{X},k)_{p+2}
$$
is not surjective, then $h^0(\cOO_\Gamma(2H-X\cap\Gamma)) > p$.
In particular, if $X$ satisfies property $\boldN_{2,p}$ and 
$\Gamma$ is a rational normal curve in $\PP^r$ not contained in $X$, 
then $$\length(X\cap \Gamma) < 2r+1-p.$$

\proof The first part follows from the conclusion of
\ref{secants} and the cohomology of a twist of 
the short exact sequence
$$
0\rTo \cII_\Gamma \rTo\cII_{\Gamma\cap X}\rTo 
\cOO_\Gamma(-X\cap\Gamma)\rTo 0.
$$
For the second part observe that the restriction map
$$
\Tor^S_{p}(I_{X\cup\Gamma},k)_{p+2}\rTo \Tor^S_{p}(I_{X},k)_{p+2}
$$
is not surjective since every element of 
$\Tor^S_{p}(I_{X\cup\Gamma},k)_{p+2}$ is represented by a syzygy
among the quadrics in $I_{X\cup\Gamma}$, which are a subset of
those in $I_\Gamma$. But a non-trivial syzygy
among the quadrics of $\Gamma$ involves all quadrics containing
$\Gamma$ and thus its syzygy variety is all of $\Gamma$ (see for instance
Ehbauer [1994] or Eisenbud-Popescu [1999]). Thus the map of 
Tor's is not surjective if $X$ is not  contained in $\Gamma$.
We may conclude now by \ref{secants} since
$h^0(\cOO_\Gamma(2H-X\cap\Gamma))=h^0(\cOO_{\PP^1}(2r-\length(X\cap\Gamma)))
> p$ if and only if $\length(X\cap \Gamma) < 2r+1-p$.\Box

\remark{binary-curve} In the special case where
both $X$ and $\Gamma$ are rational normal curves in $\PP^r$,
\ref{secants-curve} yields that $X$ and $\Gamma$ can meet at most in
$2r+1-(r-1)=r+2$ points. In case of equality, the union of the two
rational normal curves is a degeneration of a canonical curve in
$\PP^r$ (a so called ``binary'' curve). See also Eisenbud-Harris [1992],
or Diaz [1986] and Giuffrida [1988] for related results.

\remark{fail-secants-curve} The second part of \ref{secants-curve}
fails if, for instance, $\Gamma$ does not span all of $\PP^r$. Here
is an easy counterexample for $p=2$: Let $X$ be the cone
over a rational normal curve in $\PP^4$, say
$$
X=\big\{x\mid \rank\pmatrix{x_1 & x_2 &x_3\cr x_2 & x_3 & x_4}\le 1\big\}
\subset\PP^5=\PP^5(x_0,x_1,x_2,x_3,x_4,x_5)
$$
and let $\Lambda=\{x_2=x_3=x_5=0\}\subset\PP^5$.
Then $X\cap\Lambda=\{x_2=x_3=x_5=x_1x_4=0\}$ which is
a degenerate conic (union of two lines). Now if $\Gamma$
is any smooth conic in $\Lambda$, then $\length(X\cap\Gamma)=4$,
whereas \ref{secants-curve} gives an upper bound $2\cdot 2+1 -2=3$.

\medskip
\ref{secants} has numerous applications. We list in the sequel the
most interesting ones. The first one is Green's syzygetic Castelnuovo
lemma (see also Ehbauer [1994], Yanagawa [1994], and Eisenbud-Popescu
[1999]):

\corollary{castelnuovo}
Let $X\subset\PP^r$ be a finite subscheme which contains a
subscheme of length $r+3$ in linearly general position. 
If $\Tor_{r-2}(I_X, k)_r\ne 0$, then $X$ lies on a (unique)
smooth rational normal curve.

\proof Let $X'\subset X$ be a subscheme of length $r+3$ in 
linearly general position, and let $\Gamma\subset\PP^r$ 
be the unique rational normal curve containing $X'$ (see
Eisenbud-Harris [1992], or Eisenbud-Popescu [2000]).
Suppose that $X$ is not contained in $\Gamma$. Then 
the hypotheses of \ref{secants}  are 
satisfied for the scheme $X$ and the rational normal curve $\Gamma$. 
More precisely, as in the proof of \ref{secants-curve},
the restriction map
$$
\Tor^S_{r-2}(I_{X\cup\Gamma},k)_{r}\rTo \Tor^S_{r-2}(I_{X},k)_{r}
$$
is not surjective if $X$ is not  contained in $\Gamma$. We deduce from 
\ref{secants}
that $h^0(\cOO_\Gamma(2H-(X\cap\Gamma))) > r-2$.
This translates into $h^0(\cOO_{\PP^1}(2r-\length(X\cap\Gamma))) > r-2$,
whence $\length(X\cap\Gamma)\le r+2$, a contradiction since
$X\cap\Gamma$ already contains $X'$ of length $r+3$. It follows that
$X$ is contained in $\Gamma$ and this concludes the proof of
the corollary.\Box

Similarly the above techniques yield the following amusing fact (see
Eisenbud-Hulek-Popescu [2003] for more details and a better bound):

\proposition{two-veronese}
Two Veronese surfaces in $\PP^5$ whose
intersection is zero-dimensional meet in a scheme of degree at most $12$.\Box

Coble [1922] and Conner [1911] show that it is possible to realize Veronese
surfaces meeting in $10$ points, and describe such collections of points in
terms of association. A detailed analysis of the possibilities of how two
Veronese surfaces can intersect in a scheme of finite length can be found in
our paper Eisenbud-Hulek-Popescu [2003]. We show that the length of such a
scheme is at most $10$. Moreover, we prove that in the case of transversal
intersection two Veronese surfaces meet in either $10$ points (in which case
they lie on a common quadric) or in at most $8$ points. We give there also a
modern account of some of Coble and Conner's results.
                                                                                
\smallskip
Similar results hold for zero-dimensional intersections of scrolls:

\proposition{two-scrolls} Let $X$ and $\Gamma$ be two nondegenerate 
scrolls of dimensions $m$ and $n$, respectively, in $\PP^r$ with $m\le n$
and such that $X\cap\Gamma$ is a zero dimensional scheme. Then 
$\length(X\cap\Gamma)\le nr+m-{n\choose 2}+1$. 

\proof One applies \ref{secants} as in \ref{secants-curve}. The knowledge
of the number of independent quadrics in the ideal of a scroll 
(use for instance the Eagon-Northcott complex) gives then the claimed bound 
via direct computation.\Box

\references
\parindent=0pt 
\frenchspacing 
\item{} A.B.~Coble: Associated sets of points,
{\sl Trans. Amer. Math. Soc.} {\bf 24} (1922), no.1, 1--20.
\medskip

\item{} S.~Barcanescu, N.~Manolache: Betti numbers
of Segre-Veronese singularities, {\sl Rev. Roumaine
Math. Pures Appl.} {\bf 26} (1981), no. 4, 549--565.
\medskip

\item{} Th.~Bauer, T.~Szemberg: Higher order embeddings
of abelian varieties, {\sl Math. Ann.} {\bf 224} (1997), no. 3, 449-455.
\medskip

\item{} Th.~Bauer, T.~Szemberg: Primitive higher order embeddings of abelian surfaces,
{\sl Trans. Amer. Math. Soc.} {\bf 349} (1997), no. 4, 1675--1683.
\medskip

\item{} D.~Bayer, M.~Stillman: 
A criterion for detecting $m$-regularity. 
{\it Invent. Math.} {\bf 87}, (1987), no. 1, 1--11.
\medskip
 
\item{} G.~Caviglia: Bounds on the Castelnuovo-Mumford regularity
of tensor products, preprint 2003.
\medskip

\item{} A.B.~Coble: Associated sets of points,
{\sl Trans. Amer. Math. Soc.} {\bf 24} (1922), 1--20.
\medskip

\item{} J.R.~Conner: Basic systems of rational norm-curves,
{\sl Amer J. Math.} {\bf 32} (1911), no. 2, 115--176.
\medskip

\item{} S.~Diaz: Space curves that intersect often,
{\sl Pacific J. Math.} {\bf 123} (1986), no. 2, 263--267.
\medskip

\item{} I.~Dolgachev, D.~Ortland: Points sets in
projective spaces and theta functions, {\sl Ast\'erisque}
{\bf 165}, (1988).
\medskip

\item{} S.~Ehbauer: 
Syzygies of points in projective space and applications,
in {\sl ``Zero-dimensional schemes (Ravello, 1992)''}, 145--170, 
de Gruyter, Berlin, 1994. 
\medskip

\item{} L.~Ein, R.~Lazarsfeld: 
Syzygies and Koszul cohomology of smooth projective varieties of 
arbitrary dimension, 
{\sl Invent. Math.} {\bf 111} (1993), no. 1, 51--67. 
\medskip

\item{} D.~Eisenbud: 
{\sl Commutative Algebra with a View Toward Algebraic Geometry}, 
Springer, New York, 1995. 
\medskip 

\item{} D. Eisenbud: {\sl Geometry of Syzygies}, Springer, New York,
forthcoming book, 2004.
\medskip

\item{} D.~Eisenbud, S.~Goto:
Linear free resolutions and minimal multiplicity. 
{\sl J. Algebra}, {\bf 88} (1984), no. 1, 89--133.
\medskip

\item{} D.~Eisenbud, M.~Green, K.~Hulek, S.~Popescu:
The geometry of $2$-regular algebraic sets, preprint 2004.
\medskip

\item{} D.~Eisenbud, J.~Harris:
On varieties of minimal degree (a centennial account).
Algebraic geometry, Bowdoin, 1985 (Brunswick, Maine, 1985), 3--13,
{\sl Proc. Sympos. Pure Math.}, {\bf 46}, Part 1.
\medskip

\item{} D.~Eisenbud, J.~Harris:
Finite projective schemes in linearly general position.
{\sl J. Algebraic Geom.} {\bf 1} (1992), no. 1, 15--30.
\medskip

\item{} D.~Eisenbud, J.~Harris:
 An intersection bound for rank $1$ loci, with
applications to Castelnuovo and Clifford theory.
{\sl J. Algebraic Geom.} {\bf 1} (1992), no. 1, 31--59.
\medskip

\item{} D.~Eisenbud, K.~Hulek, S.~Popescu:
A note on the intersection of Veronese surfaces, in 
{\sl ``Proceedings of the NATO Advanced Workshop, held in Sinaia, Romania 2002''}, 
J.~Herzog, V.~Vuletescu eds, NATO Science series, Kluwer, 2003, 127--139.
\medskip

\item{} D.~Eisenbud, K.~Huneke, B.~Ulrich: Linearly presented ideals and a
subadditivity formula for the degrees of syzygies, preprint 2004.
\medskip

\item{} D.~Eisenbud, J.~Koh: 
Some linear syzygy conjectures. 
{\sl Adv. Math.} {\bf 90} (1991), no. 1, 47--76. 
\medskip

\item{} D.~Eisenbud, S.~Popescu:  Syzygy ideals for determinantal ideals
and the syzygetic Castelnuovo lemma. In {\sl ``Commutative Algebra, 
Algebraic Geometry, and Computational Methods, Hanoi 1996''}, 247--258, 
Springer Verlag 1999
\medskip

\item{} D.~Eisenbud, S.~Popescu: The projective geometry of the Gale Transform,
{\sl J. Algebra} {\bf 230} (2000), no. 1, 127--173.
\medskip

\item{} D.~Eisenbud, S.~Popescu, F.-O.~Schreyer, Ch.~Walter:
Exterior algebra methods for the minimal resolution conjecture,
{\sl Duke Math. J.} {\bf 112} (2002), no. 2, 379--395.
\medskip

\item{} R.~Fr\"oberg: On Stanley-Reisner rings. In
{\sl ``Topics in algebra, Part 2'' (Warsaw, 1988)}, Banach Center Publ., 
{\bf 26}, Part 2, 57--70, PWN, Warsaw, 1990.
\medskip

\item{} F.~Gallego, B.~Purnaprajna: Some results on rational
surfaces and Fano varieties, 
{\sl J. Reine Angew. Math.} {\bf 538}  (2001), 25--55.
\medskip

\item{} S.~Giuffrida:  
On the intersection of two curves in $\PP^r$,
{\sl Atti Accad. Sci. Torino Cl. Sci. Fis. Mat. Natur.} {\bf 122} (1988), 
no. 3-4, 139--143.
\medskip

\item{} D.~Grayson, M.~Stillman: 
{\sl Macaulay2}: A computer program designed to support computations
in algebraic geometry and computer algebra.  Source and object code
available from {\tt http://www.math.uiuc.edu/Macaulay2/}.
\medskip

\item{} M.~Green: Koszul cohomology and the geometry of projective varieties,
{\sl J. Differential Geom.} {\bf 19} (1984), no. 1, 125--171.
\medskip

\item{} M.~Green: The Eisenbud-Koh-Stillman conjecture on linear syzygies,
{\sl Invent. Math.} {\bf 136} (1999), no. 2, 411--418.
\medskip

\item{} M.~Green, R.~Lazarsfeld:
On the projective normality of complete linear series on an algebraic curve. 
{\sl Invent. Math.} {\bf 83} (1985), no. 1, 73--90.
\medskip 

\item{} M.~Green, R.~Lazarsfeld:
Some results on the syzygies of finite sets and algebraic curves,
{\sl Compositio Math.} {\bf 67} (1988), no. 3, 301--314.
\medskip

\item{} M.~Gross, S.~Popescu: Equations of $(1,d)$-polarized
abelian surfaces, {\it Math. Ann.} {\bf 310} (1998), no. 2, 333--377.
\medskip

\item{} J.~Herzog, T.~Hibi, X.~Zheng: Dirac's theorem on
chordal graphs and Alexander duality, preprint 2003.
\medskip

\item{} M.~Hochster: Cohen-Macaulay rings, combinatorics and simplicial
complexes, in ``Ring theory II'',  McDonald B.R., Morris, R. A. (eds),
{\it Lecture Notes in Pure and Appl. Math.}, {\bf 26}, M. Dekker 1977.
\medskip

\item{} T.~Jozefiak, P.~Pragacz, J.~Weyman:
Resolutions of determinantal varieties and tensor complexes
associated with symmetric and antisymmetric matrices,
{\sl Asterisque} {\bf 87-88} (1981), 109--189.
\medskip

\item{} A.~Lascoux: Syzygies des vari\'et\'es determinantales,
{\sl Adv. in Math.} {\bf 30} (1978), no. 3, 202--237.
\medskip

\item{} R.~Lazarsfeld: 
A sampling of vector bundle techniques in the study of linear series. In
{\sl ``Lectures on Riemann surfaces''} (Trieste, 1987), 500--559, 
World Sci. Publishing, Teaneck, NJ, 1989. 
\medskip

\item{} L.~Manivel: On the syzygies of flag manifolds,
{\sl Trans. Amer. Math. Soc.} {\bf 124} (1996), no. 8,  2293--2299.
\medskip

\item{} G.~Ottaviani, R.~Paoletti: Syzygies of Veronese embeddings,
{\sl Compositio Math.} {\bf 125}, no. 1,  (2001), 31--37.
\medskip

\item{} G.~Pareschi: Syzygies of abelian varieties.
{\sl J. Amer. Math. Soc.} {\bf 13} (2000), no. 3, 651--664.
\medskip

\item{} G.~Pareschi, M.~Popa: 
Regularity on abelian varieties I,
{\sl J. Amer. Math. Soc.} {\bf 16} (2003), no. 2, 285--302.
\medskip

\item{} P.~Pragacz, J.~Weyman: Complexes associated with trace and evaluation.
Another approach to Lascoux's resolution, {\sl Adv. in Math.} {\bf 57}
(1985), no. 2, 163--207.
\medskip

\item{} E.~Rubei: On syzygies of Segre embeddings, 
{\sl Proc. Amer. Math. Soc.}  {\bf 130}  (2002),  no. 12, 3483--3493.
\medskip

\item{} E.~Rubei: A result on resolutions of Veronese embeddings,
preprint {\tt math.AG/0309102}.
\medskip

\item{} E.~Rubei: Resolutions of Segre embeddings of projective 
spaces of any dimension, preprint {\tt math.AG/0404417}.
\medskip

\item{} N.~Shepherd-Barron, Unstable vector bundles and linear systems on
surfaces in characteristic $p$, {\sl Invent. Math.} {\bf 106} (1991),
no. 2, 243--262.
\medskip

\item{} V.~Shokurov: The Noether-Enriques theorem on canonical
curves, {\sl Mat. Sbornik} {\bf 15} (1971), 361--403 (engl. transl).
\medskip

\item{} A.J.~Sommese: Hyperplane sections of projective surfaces I. 
The adjunction mapping, {\sl Duke Math. J.} {\bf 46} (1979), no. 2, 377--401.
\medskip

\item{} A.J.~Sommese, A.~Van de Ven:
On the adjunction mapping,
{\sl Math. Ann.} {\bf 278} (1987), no. 1-4, 593--603.
\medskip

\item{} E.~Spanier: {\sl Algebraic Topology},
McGraw-Hill Book Co., 1966.
\medskip

\item{} R.~Stanley: {\sl Combinatorics and Commutative Algebra}, Second
edition, Progress in Math. {\bf 41}, Birkh\"auser, 1996.
\medskip

\item{} H.~Terakawa:
The $d$-very ampleness on a projective surface in positive characteristic,
{\sl Pacific J. of Math.} {\bf 187} (1999), no. 1, 187--198. 
\medskip

\item{} P.~Vermeire: Some results on secant varieties leading to a geometric
flip construction, {\sl Compositio Math.} {\bf 125} (2001), no. 3, 263--282.
\medskip

\item{} K.~Yanagawa: Some generalizations of Castelnuovo's lemma
on zero-dimensional schemes. {\sl J. Alg.} {\bf 170} (1994), no. 2,
429--431.
\medskip

\vskip .6cm
\widow {.1}
\vbox{\noindent Author Addresses:
\medskip
\noindent{David Eisenbud}\par
\noindent{Department of Mathematics, University of California, Berkeley,
Berkeley CA 94720}\par
\noindent{\tt de@msri.org}
\medskip
\noindent{Mark Green}\par
\noindent{Department of Mathematics, University of California at Los Angeles,
Los Angeles CA 90095-1555}\par
\noindent{\tt mlg@math.ucla.edu}
\medskip
\noindent{Klaus Hulek}\par
\noindent{Institut f\"ur Mathematik, Universit\"at Hannover, D-30060 Hannover,
Germany}\par
\noindent{\tt hulek@math.uni-hannover.de}
\medskip
\noindent{Sorin Popescu}\par
\noindent{Department of Mathematics, SUNY at Stony Brook,
Stony Brook, NY 11794-3651}\par
\noindent{\tt sorin@math.sunysb.edu}\par
}

\end